\documentclass{amsart}
\usepackage{graphicx,latexsym,amssymb}

 \newtheorem{thm}{Theorem}[section]
 \newtheorem{cor}[thm]{Corollary}
 \newtheorem{lem}[thm]{Lemma}
 \newtheorem{prop}[thm]{Proposition}
 \newtheorem{defn}[thm]{Definition}
 \newtheorem{rem}[thm]{Remark}
 \newtheorem{example}[thm]{Example}

\newcommand{\C}{\mathbb{C}}
\newcommand{\PP}{\mathbb{P}}
\newcommand{\RR}{\mathbb{R}}
\newcommand{\OO}{\mathbb{O}}
\newcommand{\s}{\sigma _}
\newcommand{\ben}{\begin{enumerate}}
\newcommand{\een}{\end{enumerate}}
\newcommand{\ble}{\begin{lem}}
\newcommand{\ele}{\end{lem}}
\newcommand{\bth}{\begin{thm}}
\renewcommand{\eth}{\end{thm}}
\newcommand{\bpr}{\begin{prop}}
\newcommand{\epr}{\end{prop}}
\newcommand{\bco}{\begin{cor}}
\newcommand{\eco}{\end{cor}}
\newcommand{\bde}{\begin{defn}}
\newcommand{\ede}{\end{defn}}
\newcommand{\brem}{\begin{rem}}
\newcommand{\erem}{\end{rem}}
\newcommand{\bexm}{\begin{example}}
\newcommand{\eexm}{\end{example}}

\newcommand{\mytable}[1]{\goodbreak\begin{table}[!ht]#1\end{table}}
\begin{document}

\title[Braid monodromy type and rational transformations]{Braid Monodromy Type and Rational Transformations of Plane Algebraic Curves}

\author[Kaplan, Shapiro and Teicher]{Shmuel Kaplan, Alexander Shapiro and Mina Teicher}

\address{Shmuel Kaplan, Alexander Shapiro and Mina Teicher,
Department of Mathematics and Statistics, Bar-Ilan University, 
Ramat-Gan 52900, Israel} \email{[kaplansh, sapial, teicher]@math.biu.ac.il} 

\thanks{This work is part of the first author Ph.D Thesis in Bar-Ilan university.}
\thanks{Partially supported by  EU-network HPRN-CT-2009-00099(EAGER) , (The Emmy Noether Research Institute for Mathematics and the Minerva Foundation of Germany), the Israel Science Foundation grant \# 8008/02-3 (Excellency Center "Group  Theoretic Methods in the Study of Algebraic Varieties").} 
\thanks{The authors wish to thank Prof. Victor Vinnikov for helpful advices.}

%%% ----------------------------------------------------------------------
\maketitle
%%% ----------------------------------------------------------------------

\begin{abstract}
We combine the newly discovered technique, which computes explicit formulas for the image of an algebraic curve under rational transformation, with techniques that enable to compute braid monodromies of such curves. We use this combination in order to study properties of the braid monodromy of the image of curves under a given rational transformation. A description of the general method is given along with full classification of the images of two intersecting lines under degree $2$ rational transformation. We also establish a connection between degree $2$ rational transformations and the local braid monodromy of the image at the intersecting point of two lines. Moreover, we present an example of two birationally isomorphic curves with the same braid monodromy type and non diffeomorphic real parts.
\end{abstract}

\section*{Introduction}

The braid monodromy is a powerful tool in the study of algebraic surfaces and curves. There exists several algorithms for computing braid monodromy for many types of curves. Usually one considers algebraic curves up to birational isomorphisms, therefore it is natural to consider the effect that a rational transformation has on the braid monodromy of a curve. Recently a new algorithm for computing the explicit image of a given algebraic curve under rational transformation was obtained. Hence, we consider the combination of these two techniques and study the braid monodromy of the image of a curve under a rational transformation. In particular, it is interesting to study the braid monodromy of an algebraic curve under a rational transformation which resolves the curve's singularities. 

In this paper we lay out the basics of the technique as follows:
In Chapter \ref{sec:Braid group preliminaries} we recall the notions and definitions of braid group, half-twists and braid monodromy. In Chapter \ref{sec:Rational transformation of the complex line} we present explicit formulas for the image of a complex line under a given rational transformation. 
We establish the connection between a rational transformation and the local braid monodromy of the image of the intersection point of two intersecting lines under this rational transformation. We present a full classification of the global braid monodromy for the image of two intersecting lines under degree $2$ rational transformation In Chapter \ref{sec:Braid monodromy of the image of two intersecting lines}. 
Chapter \ref{sec:Rational transformations of plane algebraic curve} explains a new technique which allows to find the image of curve of any degree under any rational transformation, and we give an example of degree $3$. We conclude with the computation of the global braid monodromy for the image of de-singularized curve of degree $4$.

%%%%%%%%%%%%%%%%%%%%%%%%%%%%%%%%%%%%%%%%%%%%%%%%%%%%%%%%%%%%%%%%%%%%%%%%%%%%%%%%%%%%%%%%%%%%%%%%%%%%%%%%%%%%%%%%%%%%%%%%%%%%%%
\section{Braid group preliminaries} \label{sec:Braid group preliminaries}

In this chapter we recall the definition of the braid group, some of its important elements and the braid monodromy. Readers who are interested in braid group could find more information in \cite{Artin,Birman,Dehornoy}. For information about braid monodromy we suggest readers to consult \cite{BGTI,BGTII}.

\subsection{The braid group}

\begin{defn} \label{def:Artin Braid group}
\emph{Artin's braid group $B_n$} is the group generated by $\{\sigma _1,...,\sigma _{n-1}\}$ subjected to the relations $\sigma _i \sigma _j=\sigma _j \sigma _i$ where $|i-j| \geq 2, \sigma _i \sigma _{i+1}\sigma _i=\sigma _{i+1}\sigma _i \sigma _{i+1}$ for all $i=1,...,n-2$.
\end{defn}

We distinguish some important elements in the braid group $B_n$ which are called \emph{half-twists}. Half-twists are actually the elements of the conjugacy class of any of the generators $\sigma_i$ (it is known that all generators of the braid group are conjugated to one another). 

Since the easiest way to describe and work with half-twists is based on a topological equivalent definition for the braid group, we bring it here.

$ $\\
Let $D$ be a closed disc, and $K=\{k_1,\cdots,k_n\}$ a finite set such that $K \subset int(D)$. 
\bde
Let $\mathcal{B}$ be the group of all diffeomorphisms $\beta$ of $D$ such that $\beta(K)=K$, $\beta 
|_{\partial D}={\rm Id}|_{\partial D}$. For $\beta _1,\beta _2 \in \mathcal{B}$ we say that $\beta _1$ is
\emph{equivalent} to $\beta _2$ if $\beta _1$ and $\beta _2$ induce the same automorphism of $\pi
_1 (D \setminus K,u)$, where $u$ is a point on $\partial D$. The quotient of $\mathcal{B}$ by this equivalence relation is called 
\emph{the braid group $B_n[D,K]$} ($n=|K|$). The elements of $B_n[D,K]$ are called \emph{braids}.
\ede

Now, let $D,K,u$ be as above. Let $a,b$ be two points of $K$. We denote $K_{a,b}=K \setminus
\{a,b\}$. Let $\sigma$ be a simple path in $D \setminus (\partial D \cup K_{a,b})$ connecting
$a$ with $b$. Choose a small regular neighborhood $U$ of $\sigma$ and an orientation
preserving diffeomorphism $f:\RR ^2 \to \C$ such that $f(\sigma )=[-1,1]$, $f(U)=\{z \in \C \ | \ %%@
|z| <2\}$. 

Let $\alpha (x)$, $0 \leq x$ be a real smooth monotone function such that:

$\alpha (x)=\left \{ 
\begin{array}{ll}
1, & 0 \leq x \leq \frac{3}{2} \\ 
0, & 2 \leq x
\end{array} \right.$

Define a diffeomorphism $h:\C \to \C$ as follows: for $z=re^{i\varphi} \in \C$ let
$h(z)=re^{i(\varphi +\alpha (r)\pi )}$

For the set $\{z \in \C \ | \ 2 \leq |z|\}$,  $h(z)={\rm Id}$,
and for the set $\{z \in \C \ | \ |z|\leq \frac{3}{2}\}$, $h(z)$ is a rotation by 
$180 ^{\circ}$ in the positive direction.

Considering $(f \circ h \circ f^{-1})|_D$ (we will compose from left to
right) we get a diffeomorphism of $D$ which switches $a$ and $b$ and is the identity on $D
\setminus U$. Thus it defines an element of $B_n[D,K]$.

\begin{figure}[h] 
\begin{center}
\includegraphics[scale=0.45]{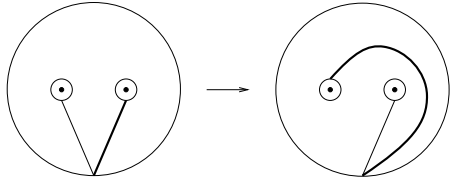}
\caption{The switch of two generators of $\pi _1(D \setminus K,u)$ induced by a half-twist.}
\label{fig:action on pi_1}
\end{center}
\end{figure}

The diffeomorphism $(f \circ h \circ f^{-1})|_D$ defined above induces an automorphism on $\pi _1(D \setminus K,u)$,
that switches the position of two generators of $\pi _1(D \setminus K,u)$, as
can be seen Figure~\ref{fig:action on pi_1}.

\bde
Let $H(\sigma)$ be the braid defined by $(f \circ h \circ f^{-1})|_D$. We call
$H(\sigma )$ the \emph{positive half-twist defined by $\sigma$}.
\ede

\begin{figure}[h] 
\begin{center}
\includegraphics[scale=0.5]{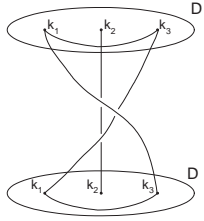}
\caption{The switch of two braid strings induced by a geometric half-twist.}
\label{fig:HT}
\end{center}
\end{figure}

The connection between the topological definition of the half-twists and the geometrical braid can be seen in Figure \ref{fig:HT}.

\subsection{The braid monodromy}

Let $C$ be a real curve in $\C^2$ of degree $n$.
Denote by $pr_1:C \to \C$ and by $pr_2:C \to \C$ the projections
to the first and second coordinate, defined in the obvious way.
For $x \in \C$ we denote $K(x)$ the projection of the points in
$C$ which lie with $x$ as their first coordinate to the second
coordinate (i.e., $K(x)=pr _2(pr _1^{-1}(x))$).

Let $N \subset \C$ be the set $N=N(C)=\{x \in\C \ | \
|K(x)|<n\}=\{x_1, \cdots ,x_p\}$. We restrict ourselves only to
the cases where $N$ is finite. Take $E$ to be a closed disc in $\C$
for which $N \subset E \setminus \partial E$. In addition take $D$
to be a closed disc in $\C$ for which $D$ contains all the points
$\{K(x) \ | \ x\in E\}$. That means that when restricted to $E$,
we have $C \subset E \times D$.

With these definitions in hand we may define the \emph{braid monodromy of a projective curve}:

\begin{defn} \label{def:braid monodromy}
Let $C$ be a projective curve of degree $n$ in $\C\PP^2$, $L$ be a generic line at infinity such that $|L \cap C|=n$, and $(x,y)$ is an affine coordinate system for $\C^2=\C\PP^2 \setminus L$ such that the projection of $C$ to the first coordinate is generic.
For $E,D,N$ defined as above, let $M \in \partial E \cap \RR$ be the base point of $\pi _1(E
\setminus N)$, and let $\sigma$ be an element of $\pi _1(E
\setminus N)$. To $\sigma$ there are $n$ lifts in $C$, each one of
them begins and ends in the points of $M \times K(M)$. Projecting
these lifts using $pr _2:C \to \C$ we get $n$ paths in
$D$ which begin and end in the points of $K(M)$. These induce a
diffeomorphism of $\pi _1(D \setminus K(M))$ which is the braid group $B_n$ as defined earlier. 
We call the homomorphism $\varphi :\pi _1(E \setminus N) \to B_n$ the \emph{braid
monodromy} of $C$ with respect to $L,E \times D,pr _1,$ and $M$.
\end{defn}

Let us fix an ordered set of generators $\left <\gamma _1, \cdots , \gamma _p \right >$ for $\pi_1(E \setminus N)$, where $p=|N|$. This set induce a $p$-tuple defined by $\left < \varphi(\gamma_1), \cdots , \varphi(\gamma _p)\right>$. We call this $p$-tuple the \emph{braid monodromy factorization} of $C$.

\begin{defn}
Let $t=(t_1,...,t_p) \in B_n^p$. We say that $s=(s_1,...,s_p) \in B_n^p$ is obtained from $t$
by the Hurwitz move $R_k$ (or $t$ is obtained from $s$ by the Hurwitz move $R_k^{-1})$ if: \\
$s_i=t_i$ for $i \neq k,k+1$ \\
$s_k=t_kt_{k+1}t_k^{-1}$ \\
$s_{k+1}=t_k$
\end{defn}

\begin{defn}
Two braid monodromy factorizations are called Hurwitz equivalent if they are obtained one from the other by a finite
sequence of Hurwitz moves and their inverses.
\end{defn}

Now we may define what the braid monodromy type of a curve is. This notion is very significant in the classification of equisingular curves as well as for the classification of surfaces.

\begin{defn}
We say that two curves are of the same \emph{BMT (Braid Monodromy Type)} if their braid monodromy factorizations are Hurwitz equivalent, up to at most one simultaneous conjugation of all elements of the first factorization by the same braid.
\end{defn}

\begin{thm} \cite{KuTe} \label{thm:KulikovTeicher}
Let $C_1$ and $C_2$ be two curves of the same BMT. Then, $C_1$ and $C_2$ are isotopic.
\end{thm}

Following from Theorem \ref{thm:KulikovTeicher} is the next corollary:

\begin{cor} \label{cor:BMT->Equivalence}
Let $C_1$ and $C_2$ be two curves, and let $F_1$ and $F_2$ be the braid monodromy factorization of $C_1$ and $C_2$ respectively. If $F_1$ and $F_2$ are Hurwitz equivalent, then the braid monodromies of $C_1$ and $C_2$ are equivalent.
\end{cor}

%%%%%%%%%%%%%%%%%%%%%%%%%%%%%%%%%%%%%%%%%%%%%%%%%%%%%%%%%%%%%%%%%%%%%%%%%%%%%%%%%%%%%%%%%%%%%%%
\section{Rational transformation of the complex line} \label{sec:Rational transformation of the complex line}
The simplest and very illustrative case is a rational
transformation of the complex line into the complex projective plane. Three
polynomials in one variable $p_0(x)$, $p_1(x)$ and $p_2(x)$ map
the complex line $\C$ into the complex projective plane $\C\PP^2$:
$$x \mapsto ( p_0(x),p_1(x),p_2(x))$$
The image of the complex line under such transformation can be
described explicitly using the notions of the Bezout matrix and
the determinantal representation of a curve. Let us recall the
definitions.

\begin{lem}[Bezout matrix] \label{lem:Bezout matrix}
For every two polynomials in one variable $p(x)$ and $q(x)$ there
exists uniquely determined $n \times n$ matrix
$B(p,q)=(b_{i,j})^n_{i,j=0}$ such that $$p(x)q(y)-q(x)p(y)=\sum _{i,j=0}^n
b_{i,j}x^i(x-y)y^j,$$\\
where $n=max\{deg(p),deg(q)\}$.
\end{lem}
This matrix is called \emph{Bezout matrix} of the polynomials $p(x)$ and
$q(x)$. For proof see \cite{SVprept}.

\begin{lem}[Determinantal representation of a curve] \label{lem:Determinantal representation of a curve}
For every homogeneous polynomial in three variables ${\Delta}(x_0,x_1,x_2)$ of degree $m$ there exist three $m \times m$ matrices $D_0$, $D_1$ and
$D_2$ such that
$${\Delta}(x_0,x_1,x_2)= \det(x_0 D_0 +x_1 D_1 + x_2 D_2)$$
\end{lem}
We will say that this is the \emph{determinantal representation} of a
curve defined by the polynomial ${\Delta}(x_0,x_1,x_2)$. For the
proof and the classification of determinantal representations of a
curve see \cite{V93}. Now we can formulate the next theorem.

\begin{thm}[Rational image of the complex line] \label{thm:Rational image of the complex line}
Let us consider three polynomials in one variable $p_0(x)$,
$p_1(x)$ and $p_2(x)$. These polynomials define the rational
transformation of the complex line $\C$ into the complex projective plane
$\C\PP^2$ by the formula $x \mapsto (p_0(x),p_1(x),p_2(x))$. The image of the complex line is
the rational curve defined by the polynomial
$$q(x_0,x_1,x_2)= \det(x_0B(p_1,p_2) +x_1 B(p_2,p_0) + x_2 B(p_0,p_1))$$
\end{thm}
For the proof see \cite{K80}.

\begin{rem}
We consider the three polynomials $p_0(x),p_1(x),p_2(x)$ to be of the same degree $n$, and if this is not the case than we consider the higher coefficients of the polynomials with degree less than $n$ to be zeroes. This is done in order to simplify the formulas, since it is equivalent to the proper formulation which uses three homogeneous polynomials of two variables that map the complex projective line $\C\PP^1$ into the complex projective plane $\C\PP^2$.
\end{rem}

%%%%%%%%%%%%%%%%%%%%%%%%%%%%%%%%%%%%%%%%%%%%%%%%%%%%%%%%%%%%%%%%%%%%%%%%%%%%%%%%%%%%%%%%%%%%%%%%%%%%%%%%%%%%%%
\section{Braid monodromy of the image of two intersecting lines} \label{sec:Braid monodromy of the image of two intersecting lines}

In this chapter we consider the image of the curve $C$ which consists of two intersecting lines defined by $x=0$ and $y=0$, under degree $2$ rational transformation $r$ into $\C\PP^2$. This implies that the image of the curve $r(C)$ consists of two intersecting conics. We consider only the generic case where the conics do not coincide and none of the conics is degenerated to a line or a point. 

\subsection{Classification of local braid monodromy}

In this section we present all the possibilities of local braid monodromy at an intersection point of two conics.

\begin{thm} \label{thm:Local braid monodromy}
The local braid monodromy of two conics at an intersection point depends only on the multiplicity of the point. More precisely, if the multiplicity of the intersection point is $n$, then the local braid monodromy is the $n$ times full-twist of two strings.
\end{thm}

\begin{proof}
There are $4$ possibilities for the multiplicity of the intersection point of two conics: $1,2,3$ and $4$. If the multiplicity is $1$ then, there exists a small neighborhood of the intersection point where the curve is the intersection of two non-singular branches (see Appendix A, Table \ref{tab:4m1}, point $x_3$). This case was studied in \cite{BGTI} and the braid monodromy was proved to be a full twists of two strings. If the multiplicity is $2$, then at the point of intersection there is a tangency of degree $1$ (see Appendix A, Table \ref{tab:2m1,1m2}, point $x_4$). This case was also previously studied in \cite{BGTII}, and the braid monodromy was proved to be a double full twists of two strings. Using the technique suggested in \cite{BGTII} this result can easily be generalized, and so in the case where the multiplicity of the intersection point is $3$ or $4$ the braid monodromy is triple or quadruple full twists of two strings. For examples see Appendix A, Table \ref{tab:1m1,1m3}, point $x_4$ and Table \ref{tab:1m4}, point $x_3$ respectively.
\end{proof}

\begin{rem}
Note that in the case of tangency with multiplicity $n$ it is easy to generalize the results and see that the local braid monodromy at the tangency point is $n$ times the full twists of two strings.
\end{rem}

Let us consider the rational transformation 
$$(x,y) \mapsto (p_0(x,y),p_1(x,y),p_2(x,y)).$$
In order to compute the local braid monodromy at the point $(p_0(x_0,y_0)$, $p_1(x_0,y_0)$, $p_2(x_0,y_0))$, 
we assume that $p_0(x_0,y_0) \neq 0$. 
We define: 
$r_1(x,y)=\frac{p_1(x,y)}{p_0(x,y)}$, $r_2(x,y)=\frac{p_2(x,y)}{p_0(x,y)}$, and recursively \\
$D_1(x)=r_2'(x,0) \cdot \frac{1}{r_1'(x,0)}$\\
$D_n(x)=D_{n-1}'(x,0) \cdot \frac{1}{r_1'(x,0)}$\\
$E_1(y)=r_2'(0,y) \cdot \frac{1}{r_1'(0,y)}$\\
$E_n(y)=D_{n-1}'(0,y) \cdot \frac{1}{r_1'(0,y)}$\\

\begin{cor} \label{cor:BM_Criteria}
Let $(p_0(x_0,y_0),p_1(x_0,y_0),p_2(x_0,y_0))$ be one of the intersection points of the two conics at the image $r(C)$. Let $i$ be the minimal index for which $D_i(x_0) \neq E_i(y_0)$. Then, the multiplicity of the intersection point is $i+1$, and thus the local braid monodromy at this intersection point is $(i+1)$ full twists of two strings.
\end{cor}

\begin{proof}
The proof follows immediately from Theorem \ref{thm:Local braid monodromy} and from the chain formula for computing derivatives. 
\end{proof}

In order to illustrate the above corollary let us consider the rational transformation 
$$(x,y) \mapsto (p_0(x,y),p_1(x,y),p_2(x,y)).$$
which is defined by the three polynomials:
$$p_0(x,y)=1+\alpha_{10}x+\alpha_{01}y+\alpha_{20}x^2+\alpha_{02}y^2$$
$$p_1(x,y)=\beta_{10}x+\beta_{01}y+\beta_{20}x^2+\beta_{02}y^2$$ 
$$p_2(x,y)=\gamma_{10}x+\gamma_{01}y+\gamma_{20}x^2+\gamma_{02}y^2$$ 
which maps the origin to the origin. According to Corollary \ref{cor:BM_Criteria} $D_1(0) \neq E_1(0)$ implies that $\gamma_{01}\beta_{10}-\gamma_{10}\beta_{01} \neq 0$, and in this case the braid monodromy at the origin is the full twist of two strings. \\
Otherwise, if $D_2(0) \neq E_2(0)$ which implies that $\beta_{10}^3(\gamma_{02}\beta_{01}-\gamma_{01}\beta_{02})+\beta_{01}^3(\gamma_{10}\beta_{20}-\gamma_{20}\beta_{10}) \neq 0$. In this case the braid monodromy at the origin is the double full twist of two strings.  \\
Otherwise, if $D_3(0) \neq E_3(0)$ which implies that 
$\beta _{10}^5(2\gamma_{02}\beta_{01}\beta_{02}+\gamma_{01}\beta_{02}\beta_{01}\alpha_{01}-2\gamma_{01}\beta_{02}^2-\gamma_{02}\alpha_{01}\beta_{01}^2) + \beta_{01}^5(2\gamma_{10}\beta_{20}^2-2\gamma_{20}\beta_{10}\beta_{20}+\gamma_{20}\alpha_{10}\beta_{10}^2-\gamma_{10}\alpha_{10}\beta_{10}\beta_{20}) \neq 0$, then the braid monodromy at the origin is the triple full twist of two strings.
Otherwise, the braid monodromy at the origin is the four times full twists of two strings. 

\subsection{Classification of braid monodromies of two intersecting conics}

In this section we give a complete classification of the braid monodromies of two intersecting conics. 
We take a projection of $\C\PP^2$ onto $\C^2$ by choosing a generic line at infinity (i.e., the line at infinity intersects the $r(C)$ at exactly $4$ points). This implies that the real part of $r(C)$ consists of two intersecting ellipses. We choose a system of coordinates for $\C^2$ in such a way that above every point of the first coordinate there is at most one singular or branch point of $r(C)$. With this construction we can compute the braid monodromy of $r(C)$ using definition \ref{def:braid monodromy}.

\begin{lem} \label{lem:HE}
The two braid monodromy factorizations:
$$F_1=\left < \s1^2,\s2,\s3\s2\s1\s2^{-1}\s3^{-1},\s3\s2^4\s3^{-1},\s3^2\s2\s3^{-2},\s3\s2^{-1}\s1\s2\s3^{-1},\s3^2 \right>$$
$$F_2=\left < \s2,\s2^{-1}\s1^2\s2,\s3\s2^{-1}\s1\s2\s3^{-1},\s1^4,\s3\s2\s1\s2^{-1}\s3^{-1},\s3^2\s2\s3^{-2},\s3^2 \right>$$
are Hurwitz equivalent.
\end{lem}

\begin{proof}
To see this activate on $F_1$ Hurwitz moves $R_1^{-1},R_5,R_4,R_3$ and $R_4$ to get $F_2$.
\end{proof}

\begin{thm} \label{thm:global braid moodromy}
Let $C$ be a curve which consists of two intersecting lines, and let $r$ be a real rational transformation of degree $2$. Then, the braid monodromy of $r(C)$ is completely defined by the number and multiplicity of it's real self intersection points. Namely,
\ben
\item
Four intersection points of multiplicity $1$, see Table \ref{tab:4m1}.
\item
Two intersection points of multiplicity $1$ and one intersection point of multiplicity $2$, see Table \ref{tab:2m1,1m2}.
\item
Two intersection points of multiplicity $2$, see Table \ref{tab:2m2}.
\item
One intersection point of multiplicity $1$ and one intersection point of multiplicity $3$, see Table \ref{tab:1m1,1m3}.
\item
One intersection point of multiplicity $4$, see Table \ref{tab:1m4}.
\item
Two intersection points of multiplicity $1$, see Table \ref{tab:2m1}.
\item
One intersection point of multiplicity $2$, see Table \ref{tab:1m2-a}.
\een
\end{thm}

\begin{proof}
For any two intersecting lines $C$ there exists a linear isomorphism between $C$ and the two intersecting lines $x=0$ and $y=0$.
Therefore, without loss of generality, we may assume that the curve $C$ consists of the two intersecting lines $x=0$ and $y=0$. 
Since $r$ is real the image of the origin is real. Hence, there are at most $2$ imaginary self intersecting points of $r(C)$. Moreover, such imaginary points must be complex conjugated. Therefore, if there exists an imaginary intersection point, its multiplicity must be $1$.

With the above assumptions all possible combinations of real self intersecting points of $r(C)$ are listed in the theorem. For combinations $1,\cdots ,6$ any two images of $C$ under rational transformations with the same type and multiplicity of intersecting points are diffeomorphic. Hence, they induce the same braid monodromy. Therefore, it is enough to compute the braid monodromy for only one example for each combination. In appendix A we give a complete description of the braid monodromy for each example. 

For combination $7$ there exists two non diffeomorphic examples. The computation of the braid monodromy for these two examples are given in Tables \ref{tab:1m2-a} and \ref{tab:1m2-b}. By Lemma \ref{lem:HE} and Corollary \ref{cor:BMT->Equivalence} these two cases yield the same BMT, hence their braid monodromies are equivalent.
\end{proof}

In Appendix $A$ we give the $8$ examples of braid monodromy computations mentioned in the proof of Theorem \ref{thm:global braid moodromy}. Each example begins by giving the polynomials which define the rational transformation $r(x,y)=(p_0(x,y),p_1(x,y),p_2(x,y))$ and the polynomial defining the image of the curve $r(C)$ (where $C$ is given by $xy=0$) under this rational transformation. Then, we give a picture of the real part of the image and a table which contains the results of the braid monodromy computation for it. Computations of the braid monodromy were performed according to the generalization of the algorithm given in \cite{BGTII}. This generalization can be found in \cite{GeneralMT}.

\begin{example}
There exists two birationally isomorphic curves of the same BMT such that their real part are not diffeomorphic.
\end{example}

\begin{proof}
See examples $7$ and $8$ in combination with Lemma \ref{lem:HE}.
\end{proof}

%%%%%%%%%%%%%%%%%%%%%%%%%%%%%%%%%%%%%%%%%%%%%%%%%%%%%%%%%%%%%%%%%%%%%%%%%%%%%%%%%%%%%%%%%%%%%%%%%%%%%%%%%%%%%%%%%
\section{Rational transformations of plane algebraic curve} \label{sec:Rational transformations of plane algebraic curve}
In this chapter we present an algorithm for computing the image of
any algebraic curve under any rational transformation. In the general
case we consider a plane real algebraic curve $C$. Let us denote
the homogeneous polynomial in three variables that defines this
curve by ${\Delta}(x_0, x_1, x_2)$. Three homogeneous polynomials in
three variables $p_0(x_0, x_1, x_2)$, $p_1(x_0,x_1, x_2)$ and $p_2(x_0, x_1, x_2)$
define the rational transformation of the plane curve $C$ by the
formula:
$$(x_0,x_1, x_2) \mapsto ( p_0(x_0,x_1,x_2),p_1(x_0,x_1,x_2),p_2(x_0,x_1,x_2)),$$ where $(x_0,x_1,x_2)$
is a point of the curve, that is ${\Delta}(x_0,x_1,x_2)=0$.

The polynomial that defines the image of the curve $C$ under such
transformation can be found using the elimination theory along an algebraic 
curve that was formulated in \cite{SVprept}. Let us recall the basic
definitions.

We will denote the degree of the polynomials $p_0(x_0,x_1,x_2)$,
$p_1(x_0,x_1,x_2)$ and $p_2(x_0,x_1,x_2)$ by $n$ and the degree of the polynomial
${\Delta}(x_0,x_1,x_2)$ that defines the curve $C$ by $m$. According to
the Lemma \ref{lem:Determinantal representation of a curve} there exists a determinantal representation of the
polynomial ${\Delta}(x_0,x_1,x_2)$, which means that there exist three $m
\times m$ matrices $D_0$, $D_1$ and $D_2$ such that
${\Delta}(x_0,x_1,x_2)= \det(x_0D_0 +x_1D_1 +x_2D_2)$. 
There is a simple way to find explicitly a determinantal representations of a polynomial
(also in more than two variables) by "lifting" it to a
noncommutative polynomial (i.e., an element of the free
associative algebra) in the same variables; see the forthcoming
work \cite{HMcCV}. This is an almost immediate corollary of the classical results of
Sch\"{u}tzenberger \cite{S61} and Fliess \cite{F74} on realization
theory for non commutative rational functions, see \cite{BR84} for a good exposition and \cite{BMGprept} for some recent progress.

We will denote by $W_n$ the space $\C^{m \frac{n(n+1)}{2}}$ as the
space of all sets of vectors $(v_{i_1 i_2})$, where each $v_{i_1
i_2} \in \C^m$ and $0 \leq i_1 + i_2 \leq n$. Let us consider a
subspace of this space
$$V_n=\{(v_{i_1 i_2}) \in W_n | D_0 v_{i_1
i_2} + D_1 v_{(i_1+1) i_2} + D_2 v_{i_1 (i_2 + 1)} = 0 \}$$

The subspace $V_n$ plays an important role in the elimination
theory along an algebraic curve and we will call $V_n$ the principal
subspace.

\begin{lem}[Generalized Bezout matrix] \label{lem:generalized Bezout matrix}
For every two homogeneous polynomials in three variables $p(x_0,x_1,x_2)$ and $q(x_0,x_1,x_2)$ of degree $n$ there exist three
$\frac{n(n+1)}{2} \times \frac{n(n+1)}{2}$ symmetric matrices
$\beta^1=(\beta^1_{i,j})$, $\beta^2=(\beta^2_{i,j})$ and
$\beta^{12}=(\beta^{12}_{i,j})$ such that $p(x_0,x_1,x_2)q(y_0,y_1,y_2)-q(x_0,x_1,x_2)p(y_0,y_1,y_2)=$
$$\sum_{|i|=|j|=n} 
\beta^1_{i,j}x^i(x_1y_0-x_0y_1)y^j +
\beta^2_{i,j}x^i(x_2y_0-x_0y_2)y^j +
\beta^{12}_{i,j}x^i(x_1y_2-x_2y_1)y^j,$$ 
where: 

$\qquad$ $i=(i_0,i_1,i_2)$, $j=(j_0,j_1,j_2)$

$\qquad$ $|i|=i_0+i_1+i_2$, $|j|=j_0+j_1+j_2$

$\qquad$ $x=(x_0,x_1,x_2)$, $y=(y_0,y_1,y_2)$

$\qquad$ $x^i=x_0^{i_0}x_1^{i_1}x_2^{i_2}$ and $y^j=y_0^{j_0}y_1^{j_1}y_2^{j_2}$.
\end{lem}

On the $m \frac{n(n+1)}{2}$-dimensional space $W_n$ let us define
a $m \frac{n(n+1)}{2} \times m \frac{n(n+1)}{2}$ matrix $B(p,q)$:
$$B(p,q)=\beta^{12} \otimes D_0 + \beta^1 \otimes D_1 + \beta^2 \otimes
D_2$$ Let us consider the restriction of $B(p,q)$ on the principal
subspace $V_n$:
$$B'(p,q)=\mathcal{P}_{V_n}B(p,q)\mathcal{P}_{V_n}$$

\begin{thm}[Rational image of a plane curve] \label{thm:Rational image of a plane curve}
Let us consider the plane real algebraic curve $C$ defined by the
polynomial $\Delta(x_0,x_1,x_2)= \det(x_0D_0 + x_1 D_1 + x_2 D_2)$ and
three homogeneous polynomials in three variables $p_0(x_0,x_1,x_2)$,
$p_1(x_0,x_1,x_2)$ and $p_2(x_0,x_1,x_2)$. These polynomials define the
rational transformation of the curve $C$ into the complex projective plane
$(x_0,x_1,x_2) \mapsto ( p_0(x_0,x_1,x_2),$ $p_1(x_0,x_1,x_2),$ $p_2(x_0,x_1,x_2))$. If the basepoints of
the transformation do not belong to the curve then the image of
the curve is defined by the polynomial
$$q(x_0,x_1,x_2)= \det(x_0B'(p_1,p_2) +x_1 B'(p_2,p_0) + x_2 B'(p_0,p_1))$$
\end{thm}

\begin{rem} If the basepoints of the transformation belong to the
curve then we have to restrict the generalized Bezout matrices
$B'(p_i,p_j)$ on a certain subspace of $V_n$ defined by the
basepoints. For details see \cite{SVprept}.
\end{rem}

\subsection{Inversion}
One of the most important rational transformations of plane
algebraic curves is the inversion. Let us consider plane real
algebraic curve $C$ defined by the degree $m$ polynomial 
${\Delta}(x_0,x_1,x_2) = \det(x_0D_0 +x_1 D_1 + x_2 D_2)$ and the
rational transformation of this curve defined by the polynomials
\begin{center}
$
\begin{array}{lcl}
p_0(x_0,x_1,x_2) & = & x_1 x_2 \\
p_1(x_0,x_1,x_2) & = & x_0 x_2 \\
p_2(x_0,x_1,x_2) & = & x_0 x_1. \\
\end{array}
$
\end{center}

We call this transformations the \emph{inversion}. The
basepoints of the inversion are the points $(0,0,1)$, $(0,1,0)$ and $(1,0,0)$. For simplicity we will assume that the
basepoints of the inversion do not belong to the curve $C$ which means that all matrices $D_0,D_1,D_2$ are non-degenerate.

In this case, $W_2=\C^{3m}$ and the principal subspace $V_2$
consists of all vectors $(v_{00}, v_{10}, v_{01})$ such that $D_0
v_{00} + D_1 v_{10} + D_2 v_{01} =0$, where $v_{00}, v_{10},
v_{01} \in \C^m$.

It is clear that $p_1(x_0,x_1,x_2)p_2(y_0,y_1,y_2)-p_2(x_0,x_1,x_2)p_1(y_0,y_1,y_2)=x_0(x_2y_1-x_1y_2)y_0=-x_0(x_1y_2-x_2y_1)y_0$.
Therefore $$B(p_1,p_2) = \beta^{12} \otimes D_0 + \beta^1 \otimes
D_1 + \beta^2 \otimes D_2 =$$
$$\left(%
\begin{array}{ccc}
  -1 & 0 & 0 \\
  0 & 0 & 0 \\
  0 & 0 & 0 \\
\end{array}%
\right) \otimes D_0 +
\left(%
\begin{array}{ccc}
  0 & 0 & 0 \\
  0 & 0 & 0 \\
  0 & 0 & 0 \\
\end{array}%
\right) \otimes D_1 +
\left(%
\begin{array}{ccc}
  0 & 0 & 0 \\
  0 & 0 & 0 \\
  0 & 0 & 0 \\
\end{array}%
\right) \otimes D_2 =$$
$$\left(%
\begin{array}{ccc}
  -D_0 & \OO & \OO \\
  \OO & \OO & \OO \\
  \OO & \OO & \OO \\
\end{array}%
\right),$$ where $\OO$ is $m \times m$ zero matrix.

For the pair $p_0$ and $p_1$ we have:
$$p_0(x_0,x_1,x_2)p_1(y_0,y_1,y_2)-p_1(x_0,x_1,x_2)p_0(y_0,y_1,y_2)=x_1x_2y_0y_2 -
x_0x_2y_1y_2 =$$ $$=x_2(x_1y_0 - x_0y_1)y_2.$$ Therefore, $$B(p_0,p_1) =
\beta^{12} \otimes D_0 + \beta^1 \otimes D_1 + \beta^2 \otimes D_2
=$$
$$\left(%
\begin{array}{ccc}
  0 & 0 & 0 \\
  0 & 0 & 0 \\
  0 & 0 & 0 \\
\end{array}%
\right) \otimes D_0 +
\left(%
\begin{array}{ccc}
  0 & 0 & 0 \\
  0 & 0 & 0 \\
  0 & 0 & 1 \\
\end{array}%
\right) \otimes D_1 +
\left(%
\begin{array}{ccc}
  0 & 0 & 0 \\
  0 & 0 & 0 \\
  0 & 0 & 0 \\
\end{array}%
\right) \otimes D_2 =$$
$$\left(%
\begin{array}{ccc}
  \OO & \OO & \OO \\
  \OO & \OO & \OO \\
  \OO & \OO & D_1 \\
\end{array}%
\right)$$

For the pair $p_2$ and $p_0$ we have:
$$p_2(x_0,x_1,x_2)p_0(y_0,y_1,y_2)-p_0(x_0,x_1,x_2)p_2(y_0,y_1,y_2)=x_0x_1y_1y_2 -
x_1x_2y_0y_1 =$$ $$=-x_1(x_2y_0 - x_0y_2)y_1.$$ Thus, $$B(p_2,p_0) = \beta^{12}
\otimes D_0 + \beta^1 \otimes D_1 + \beta^2 \otimes D_2 =$$
$$\left(%
\begin{array}{ccc}
  0 & 0 & 0 \\
  0 & 0 & 0 \\
  0 & 0 & 0 \\
\end{array}%
\right) \otimes D_0 +
\left(%
\begin{array}{ccc}
  0 & 0 & 0 \\
  0 & 0 & 0 \\
  0 & 0 & 0 \\
\end{array}%
\right) \otimes D_1 +
\left(%
\begin{array}{ccc}
  0 & 0 & 0 \\
  0 & -1 & 0 \\
  0 & 0 & 0 \\
\end{array}%
\right) \otimes D_2 =$$
$$\left(%
\begin{array}{ccc}
  \OO & \OO & \OO \\
  \OO & -D_2 & \OO \\
  \OO & \OO & \OO \\
\end{array}%
\right)$$

The next corollary follows from the Theorem \ref{thm:Rational image of a plane curve}.
\begin{cor}[Image of plane curve under the action of inversion]
Let us consider the plane real algebraic curve of $C$ degree $m$ 
defined by the polynomial $\Delta(x_0,x_1,x_2)= \det(x_0D_0 + x_1 D_1 +
x_2 D_2)$, and the rational transformation of this curve into the 
complex projective plane $(x_0,x_1,x_2) \mapsto (x_1x_2,x_0x_2,x_0x_1)$. If the points $(1,0,0)$, $(0,1,0)$ and $(0,0,1)$ do not belong to the curve $C$ then the image of the curve is defined by the polynomial
$$q(x_0,x_1,x_2)= \det \mathcal{P}_{V_2}
\left(%
\begin{array}{ccc}
  -x_0D_0 & \OO & \OO \\
  \OO & -x_1 D_2 & \OO \\
  \OO & \OO & x_2 D_1 \\
\end{array}%
\right) \mathcal{P}_{V_2},$$ where $V_2=\{(v_{00}, v_{10},
v_{01}): v_{ij} \in \C^m, D_0 v_{00} + D_1 v_{10} + D_2 v_{01}
=0\}$.
\end{cor}

\subsection{Example of the transformation of degree 3}
Generalized Bezout matrices can be found for a rational
transformation of any degree. Let us consider now rational
transformation of degree $3$ defined by the polynomials:
\begin{center}
$
\begin{array}{lcl}
p_0(x_0,x_1,x_2) & = & x_0x_1 x_2 \\
p_1(x_0,x_1,x_2) & = & x_1^3 + a x_1 x_2^2 \\
p_2(x_0,x_1,x_2) & = & x_2^3 + b x_1^2 x_2 \\
\end{array}
$
\end{center}

For the pair $p_0$ and $p_1$ one may see that\\ $p_0(x_0,x_1,x_2)p_1(y_0,y_1,y_2)-p_1(x_0,x_1,x_2)p_0(y_0,y_1,y_2)=
-x_1^2(x_1y_0-x_0y_1)y_1y_2 -x_1x_2(x_1y_0-x_0y_1)y_1^2
-ax_1x_2(x_2y_0-x_0y_2)y_1y_2 + x_1^2(x_2y_0-x_0y_2)y_1^2.$

For the pair $p_2$ and $p_0$ one may see that\\ $p_2(x_0,x_1,x_2)p_0(y_0,y_1,y_2)-p_0(x_0,x_1,x_2)p_2(y_0,y_1,y_2)=
x_2^2(x_2y_0-x_0y_2)y_1y_2 + x_1x_2(x_2y_0-x_0y_2)y_2^2
+bx_1x_2(x_1y_0-x_0y_1)y_1y_2 - x_2^2(x_1y_0-x_0y_1)y_2^2.$

For the pair $p_1$ and $p_2$ one may see that\\ $p_1(x_0,x_1,x_2)p_2(y_0,y_1,y_2)-p_2(x_0,x_1,x_2)p_1(y_0,y_1,y_2)=
x_1^2(x_1y_2-x_2y_1)y_2^2 + (1-ab)x_1x_2(x_1y_2-x_2y_1)y_1y_2 + x_2^2(x_1y_2-x_2y_1)y_1^2 +
ax_2^2(x_1y_2-x_2y_1)y_2^2 + bx_1^2(x_1y_2-x_2y_1)y_1^2.$

Therefore $$B(p_0,p_1) =
\left(%
\begin{array}{cccccc}
  \OO & \OO & \OO & \OO & \OO & \OO \\
  \OO & \OO & \OO & \OO & \OO & \OO \\
  \OO & \OO & \OO & \OO & \OO & \OO \\
  \OO & \OO & \OO &  D_2 &  -D_1 & \OO \\
  \OO & \OO & \OO & -D_1 & -aD_2 & \OO \\
  \OO & \OO & \OO & \OO & \OO & \OO \\
\end{array}%
\right)$$
$$B(p_2,p_0) =
\left(%
\begin{array}{cccccc}
  \OO & \OO & \OO & \OO & \OO & \OO \\
  \OO & \OO & \OO & \OO & \OO & \OO \\
  \OO & \OO & \OO & \OO & \OO & \OO \\
  \OO & \OO & \OO & \OO & \OO & \OO \\
  \OO & \OO & \OO & \OO & bD_1 &  D_2 \\
  \OO & \OO & \OO & \OO &  D_2 & -D_1 \\
\end{array}%
\right)$$
$$B(p_1,p_2) =
\left(%
\begin{array}{cccccc}
  \OO & \OO & \OO & \OO & \OO & \OO \\
  \OO & \OO & \OO & \OO & \OO & \OO \\
  \OO & \OO & \OO & \OO & \OO & \OO \\
  \OO & \OO & \OO & bD_0 & \OO & D_0 \\
  \OO & \OO & \OO & \OO & (1-ab)D_0 & \OO \\
  \OO & \OO & \OO & D_0 & \OO & aD_0 \\
\end{array}%
\right)$$
and the image of the curve under this transformation can
be found from Theorem \ref{thm:Rational image of a plane curve}.

%%%%%%%%%%%%%%%%%%%%%%%%%%%%%%%%%%%%%%%%%%%%%%%%%%%%%%%%%%%%%%%%%%%%%%%%%%
\subsection{Example of braid monodromy of de-singularized curve}
The technique described in this chapter allows to compute explicit formulas for the image of an algebraic curve of any degree under any rational transformation. In particular, it seems to be interesting to study the connection between the braid monodromies of singular algebraic curves and their images under the action of de-singularizing rational transformations.

Let us consider the degree $4$ singular curve $C$ defined by the polynomial:
$$q(x,y)=x^3y-xy^3+2x^3-y^3$$

\begin{figure}[h]
\begin{center}
\includegraphics[scale=0.5,angle=0,draft=false,angle=-90]{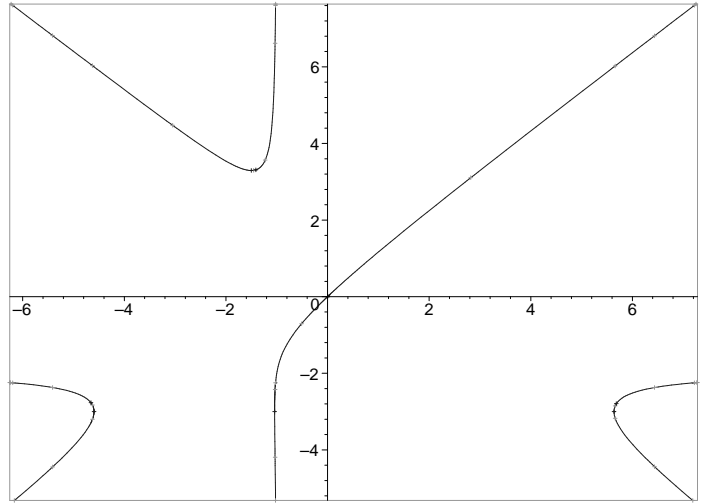}
\caption{The real part of $C$}
\label{fig:ComplicatedNode}
\end{center}
\end{figure}

Figure \ref{fig:ComplicatedNode} is the graph of the real part of $C$.

By the technique mentioned in this chapter, it is possible to find a determinantal representation for this curve. One of these determinantal representations is:

$$det \left ( \left(
\begin{array}{cccc}
2 & 0 & -2 & 0\\ 
0 & -2 & -1 & 0\\
-2 & -1 & 0 & 0\\
0 & 0 & 0 & 0
\end{array} \right ) 
+
x \left (
\begin{array}{cccc}
0 & 0 & 0 & 0\\
0 & 0 & -2 & -1\\
0 & -2 & -1 & 0\\
0 & -1 & 0 & 1
\end{array} \right ) 
+
y \left (
\begin{array}{cccc}
4 & 0 & 0 & 2\\
0 & 0 & 0 & 0\\
0 & 0 & 0 & 0\\
2 & 0 & 0 & 1
\end{array} \right ) 
\right ) $$

It is clear that some singularity of this curve occurs at infinity. This implies that it might turn to be a very complicated task to compute the braid monodromy of this curve. On the other hand, we may de-singularize this curve using inversion. The image will then be the curve $C_1$ defined by the polynomial $$q_1(x,y)=det\left ( \left(
\begin{array}{ccc}
-2 & 0 & 2\\
0 & 2 & 1\\
2 & 1 & 0
\end{array} \right ) 
+
x \left (
\begin{array}{ccc}
0 & 2 & 1\\
2 & 1 & 0\\
1 & 0 & -1
\end{array} \right ) 
+
y \left (
\begin{array}{ccc}
-1 & 0 & 1\\
0 & 0 & 0\\
1 & 0 & -1
\end{array} \right ) 
\right )=$$
$$=x^3-2y^3+x^2-y^2$$

\begin{figure}[h]
\begin{center}
\includegraphics[scale=0.5,angle=0,draft=false,angle=-90]{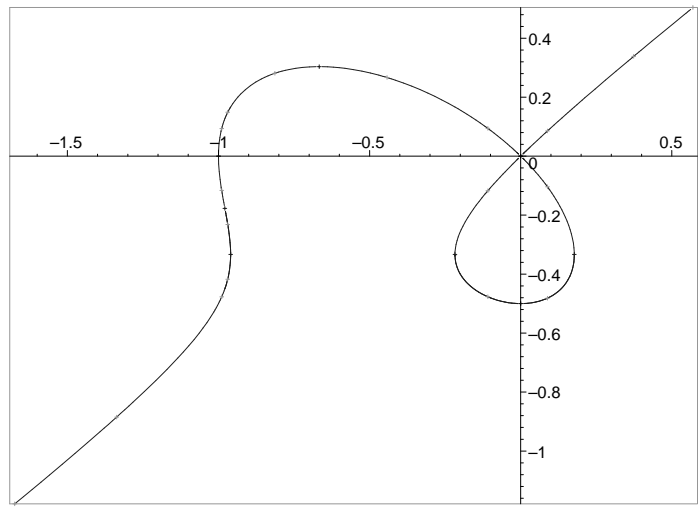}
\caption{The real part of $C_1$}
\label{fig:Node}
\end{center}
\end{figure}

It is clear that $C$ and $C_1$ are birationally isomorphic, and that $C_1$ is an almost real curve (i.e., it is defined with real coefficients and all its singular and branch points all have different real coordinates). Therefore, it is possible to compute the braid monodromy of $C_1$ using the algorithm given in \cite{BGTII}. Figure \ref{fig:Node} is the graph of the real part of $C_1$, and Table \ref{tab:NodeBM} describes its braid monodromy.

\mytable{
\begin{center}
\begin{tabular}{||l||c|c|c|c|c||}
\hline
\hline
& & & & &\\
Singular point & $x_1$ &   $x_2$ & $x_3$ & $x_4$ & $x_5$\\ 
& & & & &\\
\hline
& & & & &\\
Braid monodromy &
$\s2^{-1}\s1\s2$ & 
$\s1^2$ & 
$\s2$ &
$\s2$ &
$\s2^{-1}\s1\s2$ \\
& & & & &\\
\hline
%& & & & \\
\includegraphics[scale=1]{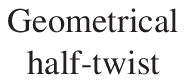} &
\includegraphics[scale=0.5]{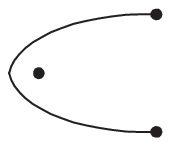} &
\includegraphics[scale=0.5]{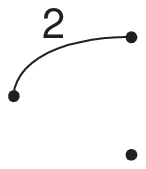} &
\includegraphics[scale=0.5]{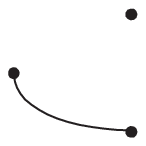} &
\includegraphics[scale=0.5]{NodeX3.eps} &
\includegraphics[scale=0.5]{NodeX1.eps} \\
\hline
\hline
\end{tabular}
\caption{Braid monodromy results for the curve $C_1$.}\label{tab:NodeBM}
\end{center}}

We saw that there are some connections between the braid monodromy of plane algebraic curves and the local braid monodromy of its image under degree $2$ rational transformation. Generalizing this connection for curves and rational transformations of higher degree is the goal of our next research.

%%%%%%%%%%%%%%%%%%%%%%%%%%%%%%%%%%%%%%%%%%%%%%%%%%%%%%%%%%%%%%%%%%%%%%%%%%%%%%%%%%%%%%%%%%%%%%%
\clearpage

\section{Appendix A - Braid monodromies for the proof of Theorem \ref{thm:global braid moodromy}}

%%%%%%%%%%%%%%%%%%%%%%%%%%%%%%%%%%%%%%%%%%%%%%%%%%%%%%%%%%%%%%%%%%%%%%%%%%%%%%%%%%%%%%%%%%%%%%%
%                                                                                             %
% Four intersection points of multiplicity 1                                                  %
%                                                                                             %
%%%%%%%%%%%%%%%%%%%%%%%%%%%%%%%%%%%%%%%%%%%%%%%%%%%%%%%%%%%%%%%%%%%%%%%%%%%%%%%%%%%%%%%%%%%%%%%

\subsection{Example 1. Four intersection points of multiplicity 1} $ $\\
$p_0(x,y)=1+x^2+y^2$, $p_1(x,y)=2x+4y$, $p_2(x,y)=2x^2+3x+y^2$ \\
$r(C)$ is defined by: $(x^2+16y^2-16y)(13x^2-12xy+4y^2+12x-8y)=0$

\begin{figure}[h] 
\begin{center}
\includegraphics[scale=0.45,angle=-90]{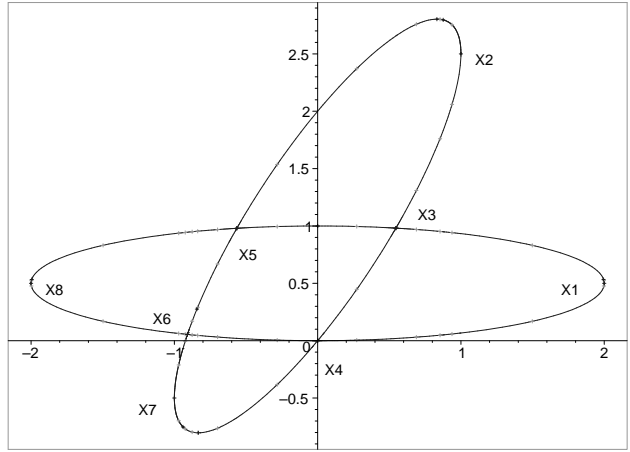}
\caption{Real part of $r(C)$ for Example $1$}
\label{fig:4m1}
\end{center}
\end{figure}
\mytable{
\begin{center}
\begin{tabular}{||l||c|c|c|c||}
\hline
\hline
& & & & \\
Singular point & $x_1$ &   $x_2$ & $x_3$ & $x_4$\\ 
& & & & \\
\hline
& & & & \\
Braid monodromy &
$\s2$ & 
$\s3\s2\s1\s2^{-1}\s3^{-1} $ & 
$\s3\s2^2\s3^{-1}$ &
$\s3^2$ \\
& & & & \\
\hline
%& & & & \\
\includegraphics[scale=1]{TitleGeometricHT.eps} &
\includegraphics[scale=0.5]{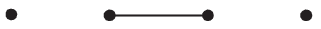} &
\includegraphics[scale=0.5]{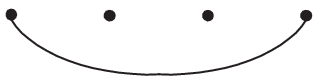} &
\includegraphics[scale=0.5]{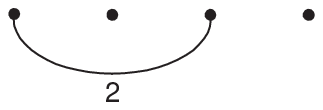} &
\includegraphics[scale=0.5]{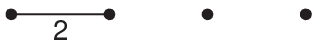} \\
\hline
\hline
& & & & \\
Singular point & $x_5$ & $x_6$ & $x_7$ & $x_8$ \\ 
& & & &\\
\hline
& & & &\\
Braid monodromy &
$\s1^2$ & 
$\s2\s1^2\s2^{-1}$  & 
$\s3\s2\s1\s2^{-1}\s3^{-1}$ &
$\s2$ \\
& & & & \\
\hline
%& & & & \\
\includegraphics[scale=1]{TitleGeometricHT.eps} &
\includegraphics[scale=0.5]{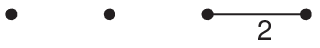} &
\includegraphics[scale=0.5]{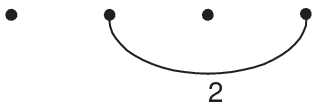} & 
\includegraphics[scale=0.5]{1x2.eps} &
\includegraphics[scale=0.5]{1x1.eps} \\
\hline
\hline
\end{tabular}
\caption{Braid monodromy results for Example $1$.}\label{tab:4m1}
\end{center}}
\newpage

%%%%%%%%%%%%%%%%%%%%%%%%%%%%%%%%%%%%%%%%%%%%%%%%%%%%%%%%%%%%%%%%%%%%%%%%%%%%%%%%%%%%%%%%%%%%%%%
%                                                                                             %
% Two intersection points of multiplicity $1$ and one of multiplicity $2$                     %
%                                                                                             %
%%%%%%%%%%%%%%%%%%%%%%%%%%%%%%%%%%%%%%%%%%%%%%%%%%%%%%%%%%%%%%%%%%%%%%%%%%%%%%%%%%%%%%%%%%%%%%%
\subsection{Example 2. Two intersection points of multiplicity 1 and one of multiplicity 2} $ $\\
$p_0(x,y)=1+x^2+y^2$, $p_1(x,y)=1.5x+2y$, $p_2(x,y)=3x^2+y^2$ \\
$r(C)$ is defined by: $(x^2+4y^2-4y)(36x^2+9y^2-27y)=0$

\begin{figure}[h] 
\begin{center}
\includegraphics[scale=0.45,angle=-90]{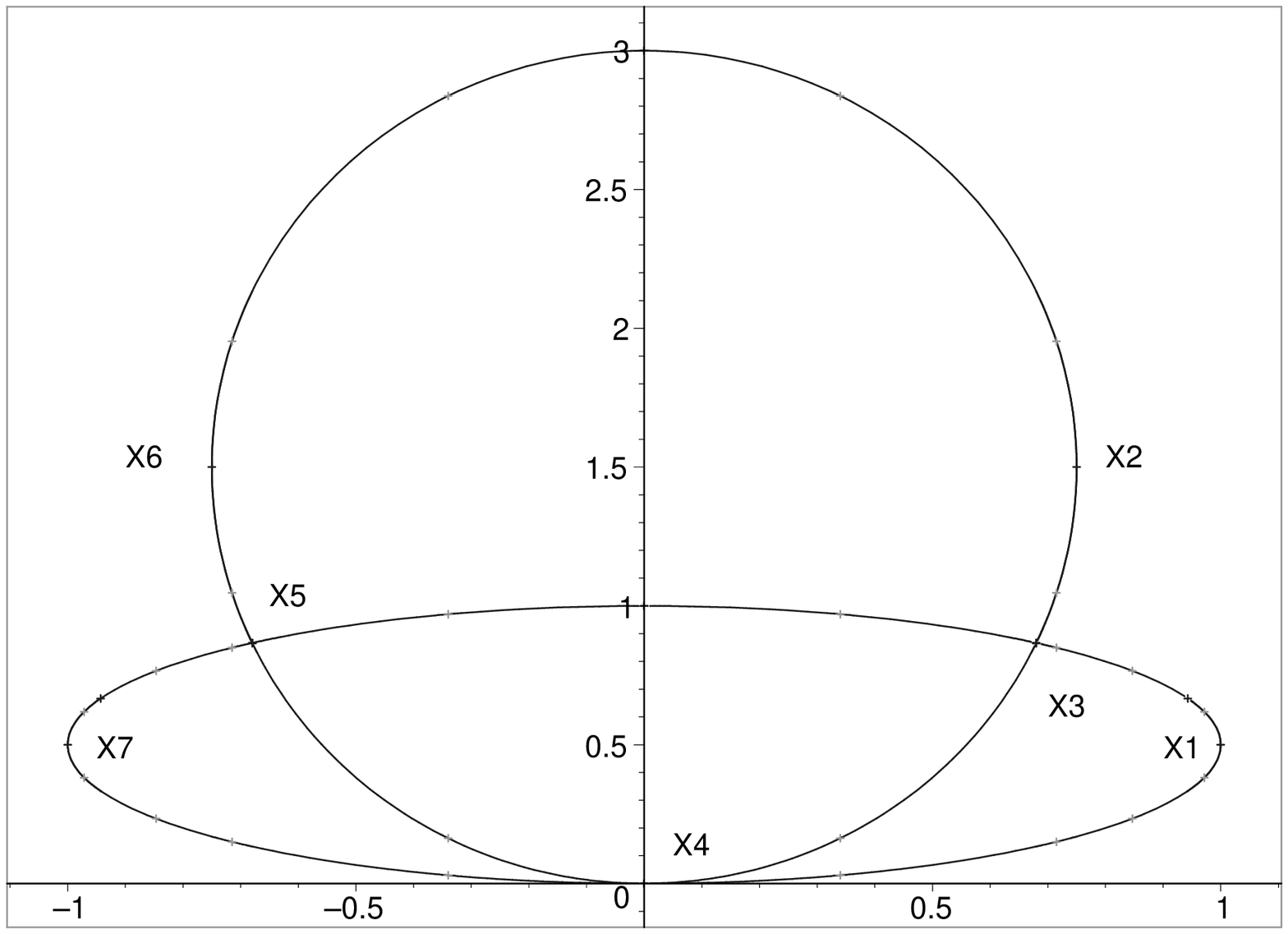}
\caption{Real part of $r(C)$ for Example $2$}
\label{fig:2m1,1m2}
\end{center}
\end{figure}

\mytable{
\begin{center}
\begin{tabular}{||l||c|c|c|c||}
\hline
\hline
& & & & \\
Singular point & $x_1$ &   $x_2$ & $x_3$ & $x_4$\\ 
& & & & \\
\hline
& & & & \\
Braid monodromy &
$\s2$ & 
$\s3\s2\s1\s2^{-1}\s3^{-1} $ & 
$\s3\s2^2\s3^{-1}$ &
$\s3^4$ \\
& & & & \\
\hline
%& & & & \\
\includegraphics[scale=1]{TitleGeometricHT.eps} &
\includegraphics[scale=0.5]{1x1.eps} &
\includegraphics[scale=0.5]{1x2.eps} &
\includegraphics[scale=0.5]{1x3.eps} &
\includegraphics[scale=0.5]{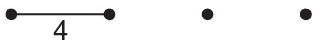} \\
\hline
\hline
& & & & \\
Singular point & $x_5$ & $x_6$ & $x_7$ & \\ 
& & & &\\
\cline{1-4}
& & & &\\
Braid monodromy &
$\s3^{-1}\s2^2\s3$ & 
$\s3^{-1}\s2^{-1}\s1\s2\s3$  & 
$\s2$ & \\
& & & & \\
\cline{1-4}
%& & & & \\
\includegraphics[scale=1]{TitleGeometricHT.eps} &
\includegraphics[scale=0.5]{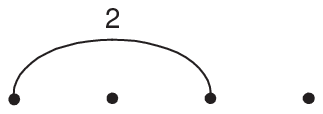} &
\includegraphics[scale=0.5]{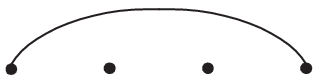} & 
\includegraphics[scale=0.5]{1x1.eps} &
\\
\hline
\hline
\end{tabular}
\caption{Braid monodromy results for Example $2$}\label{tab:2m1,1m2}
\end{center}}
\newpage

%%%%%%%%%%%%%%%%%%%%%%%%%%%%%%%%%%%%%%%%%%%%%%%%%%%%%%%%%%%%%%%%%%%%%%%%%%%%%%%%%%%%%%%%%%%%%%%
%                                                                                             %
% Two intersection points of multiplicity $2$                                                 %
%                                                                                             %
%%%%%%%%%%%%%%%%%%%%%%%%%%%%%%%%%%%%%%%%%%%%%%%%%%%%%%%%%%%%%%%%%%%%%%%%%%%%%%%%%%%%%%%%%%%%%%%
\subsection{Example 3. Two intersection points of multiplicity 2} $ $\\
$p_0(x,y)=1+x^2+y^2$, $p_1(x,y)=2x+y+x^2+y^2$, $p_2(x,y)=2x+y-x^2-y^2$ \\
$r(C)$ is defined by: $(2x^2+2y^2-2x+2y)(5x^2-6xy+5y^2-8x+8y)=0$

\begin{figure}[h] 
\begin{center}
\includegraphics[scale=0.45,angle=-90]{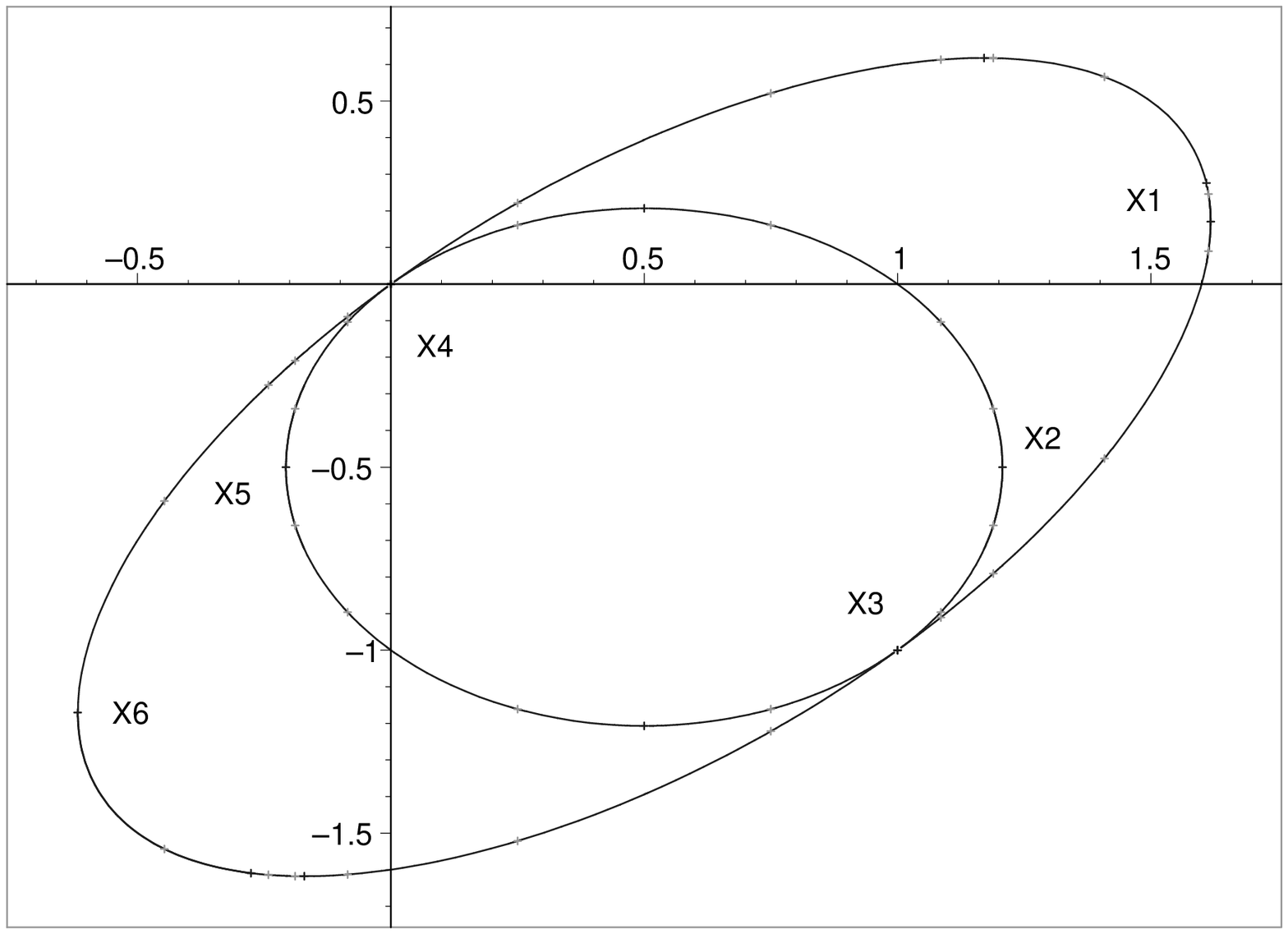}
\caption{Real part of $r(C)$ for Example $3$}
\label{fig:2m2}
\end{center}
\end{figure}

\mytable{
\begin{center}
\begin{tabular}{||l||c|c|c|c||}
\hline
\hline
& & & & \\
Singular point & $x_1$ &   $x_2$ & $x_3$ & $x_4$\\ 
& & & & \\
\hline
& & & & \\
Braid monodromy &
$\s2$ & 
$\s3\s2^{-1}\s1\s2\s3^{-1}$ & 
$\s3^4$ &
$\s1^4$ \\
& & & & \\
\hline
%& & & & \\
\includegraphics[scale=1]{TitleGeometricHT.eps} &
\includegraphics[scale=0.5]{1x1.eps} &
\includegraphics[scale=0.5]{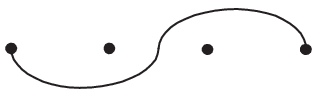} &
\includegraphics[scale=0.5]{2x4.eps} &
\includegraphics[scale=0.5]{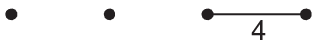} \\
\hline
\hline
& & & & \\
Singular point & $x_5$ & $x_6$ & & \\ 
& & & &\\
\cline{1-3}
& & & &\\
Braid monodromy &
$\s3^{-1}\s2\s1\s2^{-1}\s3$ & 
$\s2$ & & \\
& & & & \\
\cline{1-3}
%& & & & \\
\includegraphics[scale=1]{TitleGeometricHT.eps} &
\includegraphics[scale=0.5]{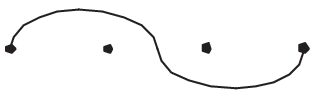} &
\includegraphics[scale=0.5]{1x1.eps} & &
\\
\hline
\hline
\end{tabular}
\caption{Braid monodromy results for Example $3$.}\label{tab:2m2}
\end{center}}
\newpage

%%%%%%%%%%%%%%%%%%%%%%%%%%%%%%%%%%%%%%%%%%%%%%%%%%%%%%%%%%%%%%%%%%%%%%%%%%%%%%%%%%%%%%%%%%%%%%%
%                                                                                             %
% One intersection point of multiplicity $1$ and one intersection point of multiplicity $3$   %
%                                                                                             %
%%%%%%%%%%%%%%%%%%%%%%%%%%%%%%%%%%%%%%%%%%%%%%%%%%%%%%%%%%%%%%%%%%%%%%%%%%%%%%%%%%%%%%%%%%%%%%%
\subsection{Example 4. One intersection point of multiplicity 1 and one intersection point of multiplicity 3} $ $\\
$p_0(x,y)=1+x^2+y^2$, $p_1(x,y)=x+x^2+y$, $p_2(x,y)=x+y-y^2$ \\
$r(C)$ is defined by: $(2x^2-2xy+y^2-x+y)(x^2-2xy+2y^2-x+y)=0$

\begin{figure}[h] 
\begin{center}
\includegraphics[scale=0.45,angle=-90]{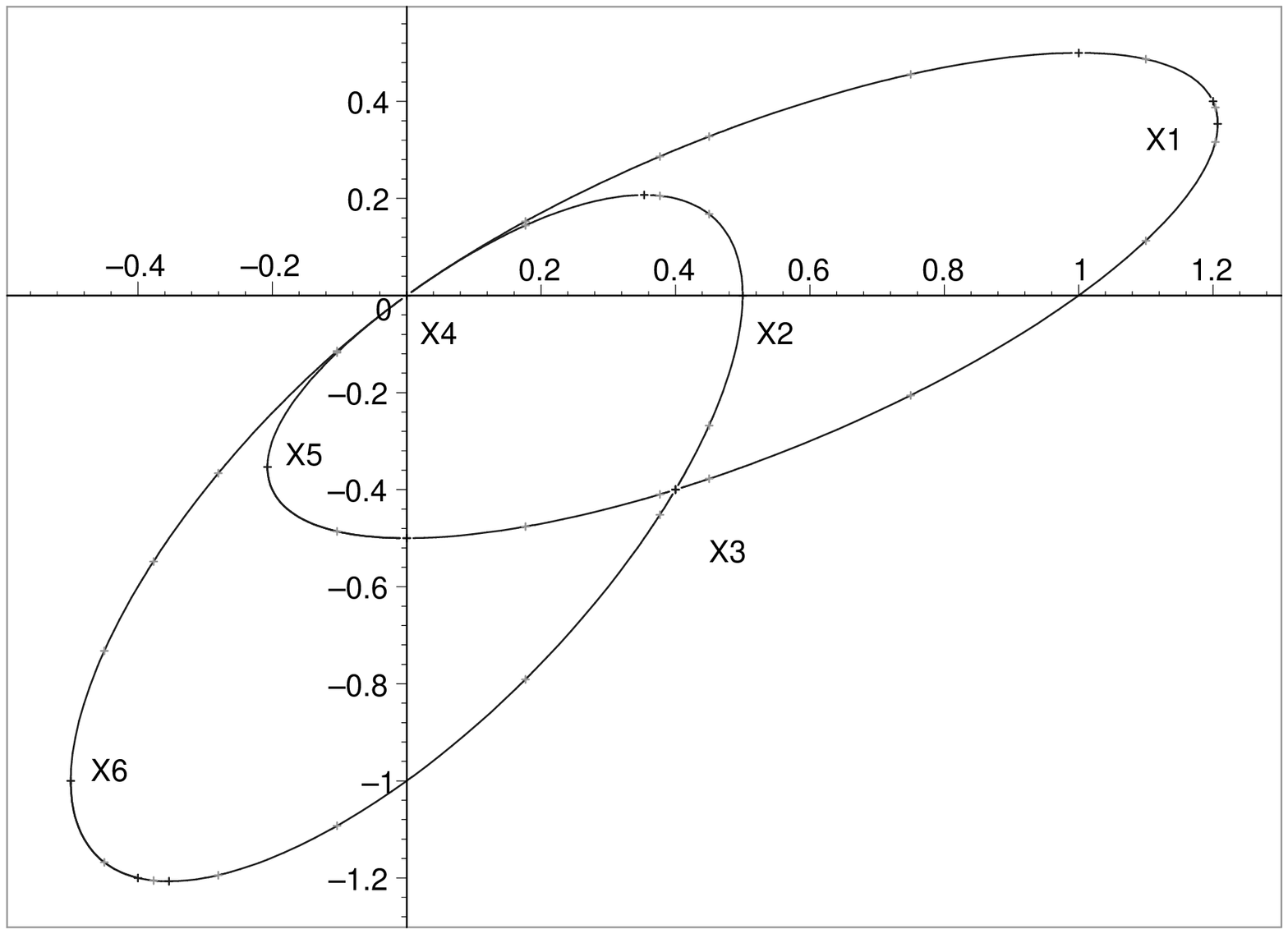}
\caption{Real part of $r(C)$ for Example $4$}
\label{fig:1m1,1m3}
\end{center}
\end{figure}

\mytable{
\begin{center}
\begin{tabular}{||l||c|c|c|c||}
\hline
\hline
& & & & \\
Singular point & $x_1$ &   $x_2$ & $x_3$ & $x_4$\\ 
& & & & \\
\hline
& & & & \\
Braid monodromy &
$\s2$ & 
$\s3\s2^{-1}\s1\s2\s3^{-1}$ & 
$\s3^2$ &
$\s1^6$ \\
& & & & \\
\hline
%& & & & \\
\includegraphics[scale=1]{TitleGeometricHT.eps} &
\includegraphics[scale=0.5]{1x1.eps} &
\includegraphics[scale=0.5]{7bX3.eps} &
\includegraphics[scale=0.5]{1x4.eps} &
\includegraphics[scale=0.5]{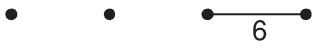} \\
\hline
\hline
& & & & \\
Singular point & $x_5$ & $x_6$ & & \\ 
& & & &\\
\cline{1-3}
& & & &\\
Braid monodromy &
$\s1^{-2}\s2\s1^2$ & 
$\s3\s2\s1\s2^{-1}\s3^{-1}$ & & \\
& & & & \\
\cline{1-3}
%& & & & \\
\includegraphics[scale=1]{TitleGeometricHT.eps} &
\includegraphics[scale=0.5]{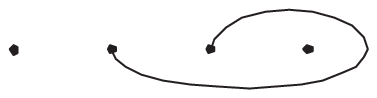} &
\includegraphics[scale=0.5]{1x2.eps} & &
\\
\hline
\hline
\end{tabular}
\caption{Braid monodromy results for Example $4$.}\label{tab:1m1,1m3}
\end{center}}
\newpage

%%%%%%%%%%%%%%%%%%%%%%%%%%%%%%%%%%%%%%%%%%%%%%%%%%%%%%%%%%%%%%%%%%%%%%%%%%%%%%%%%%%%%%%%%%%%%%%
%                                                                                             %
% One intersection point of multiplicity $4$                                                  %
%                                                                                             %
%%%%%%%%%%%%%%%%%%%%%%%%%%%%%%%%%%%%%%%%%%%%%%%%%%%%%%%%%%%%%%%%%%%%%%%%%%%%%%%%%%%%%%%%%%%%%%%
\subsection{Example 5. One intersection point of multiplicity 4} $ $\\
$p_0(x,y)=1+x^2+y^2$, $p_1(x,y)=2x+y$, $p_2(x,y)=4x^2+y^2$ \\
$r(C)$ is defined by: $(x^2+y^2-y)(4x^2+y^2-4y)=0$

\begin{figure}[h] 
\begin{center}
\includegraphics[scale=0.45,angle=-90]{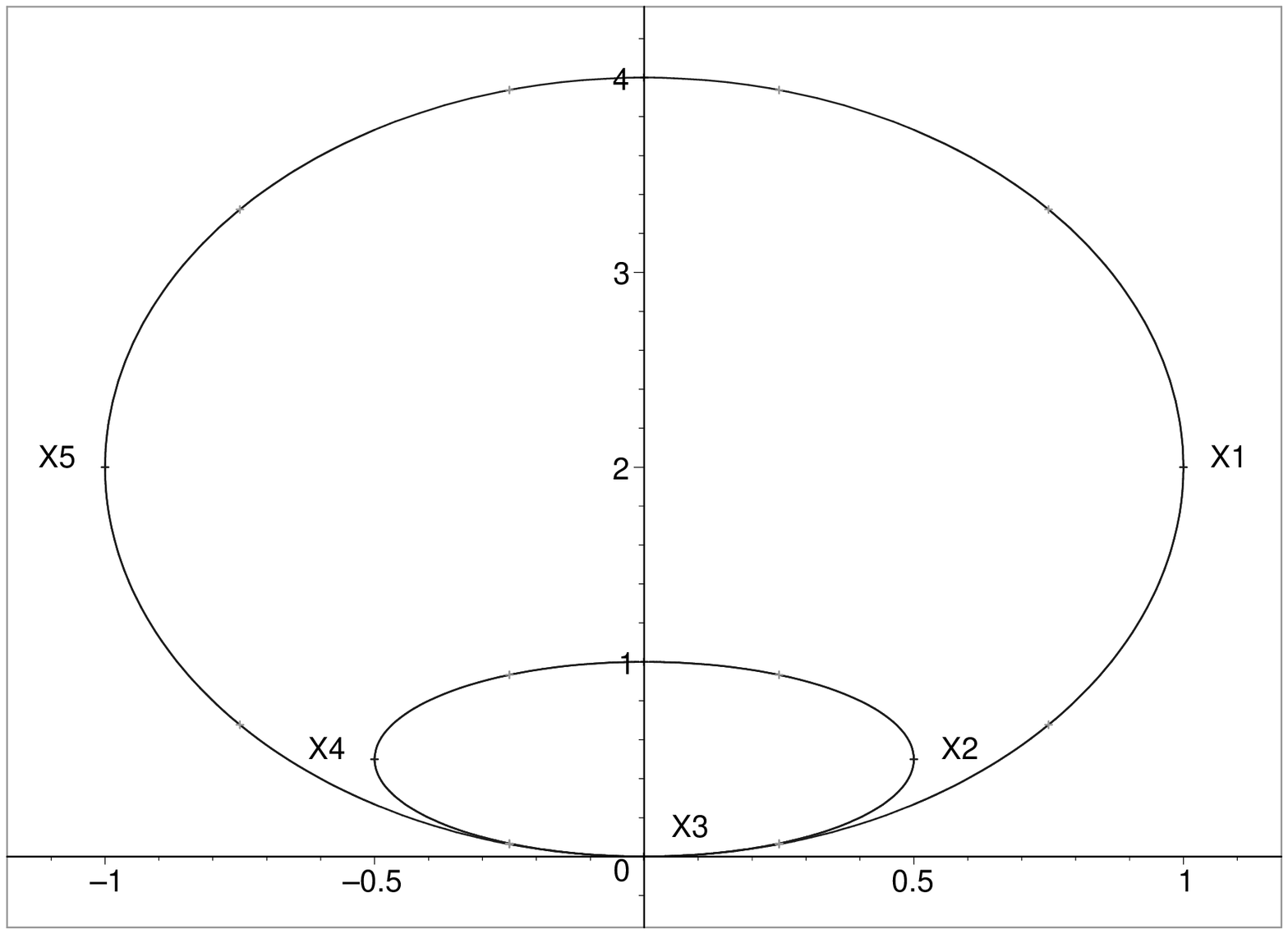}
\caption{Real part of $r(C)$ for Example $5$}
\label{fig:1m4}
\end{center}
\end{figure}

\mytable{
\begin{center}
\begin{tabular}{||l||c|c|c|c||}
\hline
\hline
& & & & \\
Singular point & $x_1$ &   $x_2$ & $x_3$ & $x_4$\\ 
& & & & \\
\hline
& & & & \\
Braid monodromy &
$\s2$ & 
$\s3\s2^{-1}\s1\s2\s3^{-1}$ & 
$\s3^8$ &
$\s3^{-3}\s2^{-1}\s1\s2\s3^3$ \\
& & & & \\
\hline
%& & & & \\
\includegraphics[scale=1]{TitleGeometricHT.eps} &
\includegraphics[scale=0.5]{1x1.eps} &
\includegraphics[scale=0.5]{7bX3.eps} &
\includegraphics[scale=0.5]{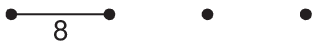} &
\includegraphics[scale=0.5]{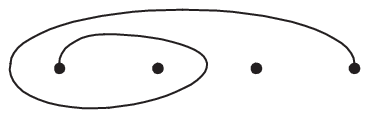} \\
\hline
\hline
& & & & \\
Singular point & $x_5$ & & & \\ 
& & & &\\
\cline{1-2}
& & & &\\
Braid monodromy &
$\s1^2\s3^{-1}\s2\s3\s2^{-1}\s3\s1^{-2}$ & 
& & \\
& & & & \\
\cline{1-2}
%& & & & \\
\includegraphics[scale=1]{TitleGeometricHT.eps} &
\includegraphics[scale=0.5]{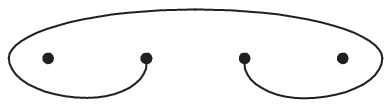} &
& & \\
\hline
\hline
\end{tabular}
\caption{Braid monodromy results for Example $5$.}\label{tab:1m4}
\end{center}}
\newpage

%%%%%%%%%%%%%%%%%%%%%%%%%%%%%%%%%%%%%%%%%%%%%%%%%%%%%%%%%%%%%%%%%%%%%%%%%%%%%%%%%%%%%%%%%%%%%%%
%                                                                                             %
% Two intersection points of multiplicity $1$                                                 %
%                                                                                             %
%%%%%%%%%%%%%%%%%%%%%%%%%%%%%%%%%%%%%%%%%%%%%%%%%%%%%%%%%%%%%%%%%%%%%%%%%%%%%%%%%%%%%%%%%%%%%%%
\subsection{Example 6. Two intersection points of multiplicity 1} $ $\\
$p_0(x,y)=1+x^2+y^2$, $p_1(x,y)=x+x^2+y$, $p_2(x,y)=x+y^2$ \\
$r(C)$ is defined by: $(x^2+y^2-y)(x^2-2xy+2y^2-x+y)=0$

\begin{figure}[h] 
\begin{center}
\includegraphics[scale=0.45,angle=-90]{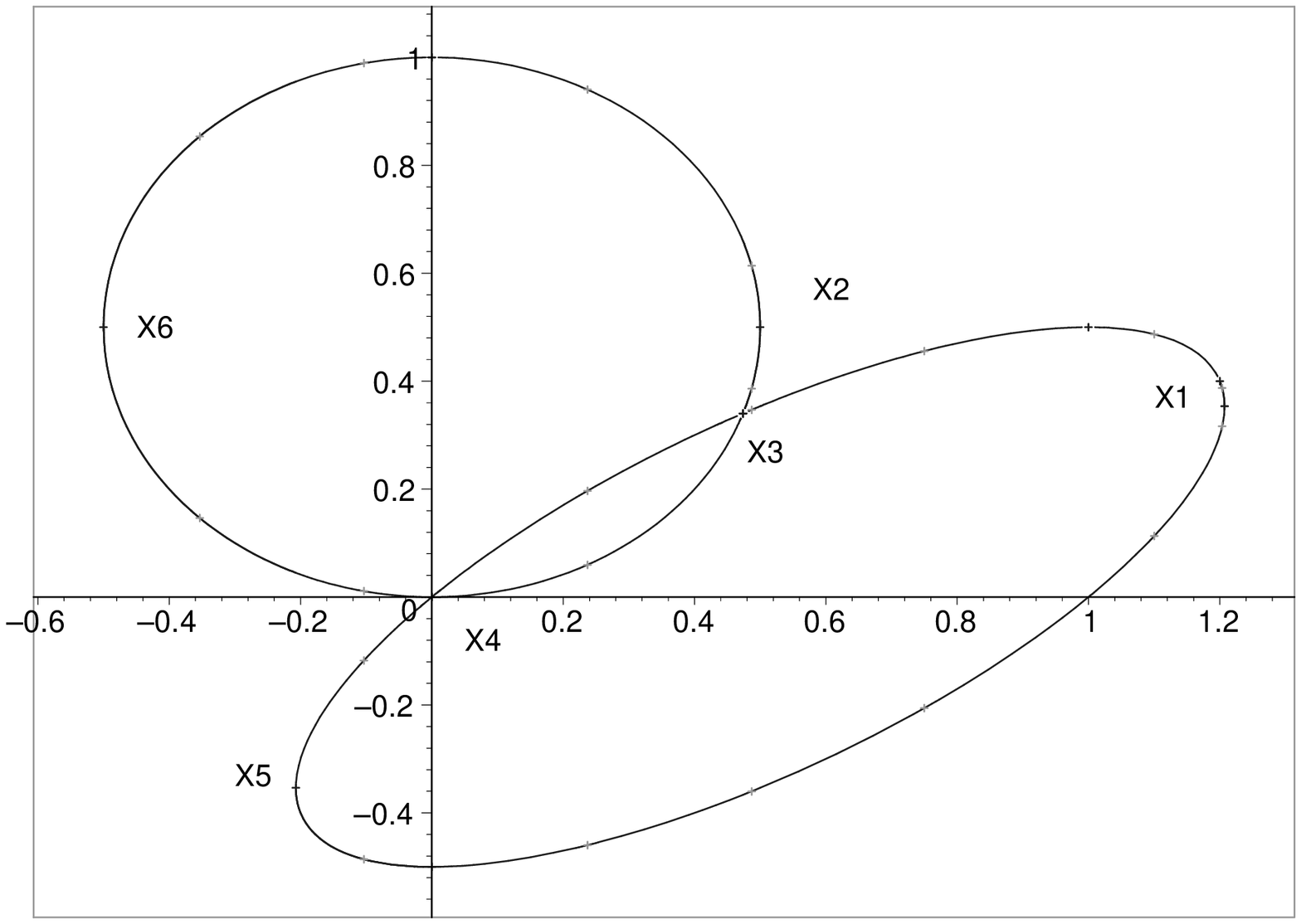}
\caption{Real part of $r(C)$ for Example $6$}
\label{fig:2m1}
\end{center}
\end{figure}

\mytable{
\begin{center}
\begin{tabular}{||l||c|c|c|c||}
\hline
\hline
& & & & \\
Singular point & $x_1$ &   $x_2$ & $x_3$ & $x_4$\\ 
& & & & \\
\hline
& & & & \\
Braid monodromy &
$\s2$ & 
$\s2^{-1}\s1^2\s2$ & 
$\s3\s2\s1\s2^{-1}\s3^{-1}$ &
$\s3\s2^2\s3^{-1}$ \\
& & & & \\
\hline
%& & & & \\
\includegraphics[scale=1]{TitleGeometricHT.eps} &
\includegraphics[scale=0.5]{1x1.eps} &
\includegraphics[scale=0.5]{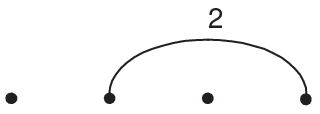} &
\includegraphics[scale=0.5]{1x2.eps} &
\includegraphics[scale=0.5]{1x3.eps} \\
\hline
\hline
& & & & \\
Singular point & $x_5$ & $x_6$ & $x_7$ & $x_8$ \\ 
& & & &\\
\hline
& & & &\\
Braid monodromy &
$\s3\s2^2\s3^{-1}$ & 
$\s3^2\s2\s3^{-2}$  & 
$\s3\s2^{-1}\s1\s2\s3^{-1}$ &
$\s3^2$ \\
& & & & \\
\hline
%& & & & \\
\includegraphics[scale=1]{TitleGeometricHT.eps} &
\includegraphics[scale=0.5]{1x3.eps} &
\includegraphics[scale=0.5]{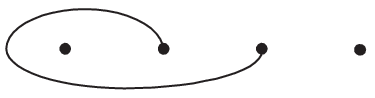} & 
\includegraphics[scale=0.5]{7bX3.eps} &
\includegraphics[scale=0.5]{1x4.eps} \\
\hline
\hline
\end{tabular}
\caption{Braid monodromy results for Example $6$}\label{tab:2m1}
\end{center}}
\newpage

%%%%%%%%%%%%%%%%%%%%%%%%%%%%%%%%%%%%%%%%%%%%%%%%%%%%%%%%%%%%%%%%%%%%%%%%%%%%%%%%%%%%%%%%%%%%%%%
%                                                                                             %
% One intersection points of multiplicity $2$ - type a                                        %
%                                                                                             %
%%%%%%%%%%%%%%%%%%%%%%%%%%%%%%%%%%%%%%%%%%%%%%%%%%%%%%%%%%%%%%%%%%%%%%%%%%%%%%%%%%%%%%%%%%%%%%%
\subsection{Example 7. One intersection points of multiplicity 2 - type a} $ $\\
$p_0(x,y)=1+x^2+y^2$, $p_1(x,y)=3x+y+2x^2-y^2$, $p_2(x,y)=3x+y-2x^2+y^2$ \\
$r(C)$ is defined by: $(x^2+y^2+x-y)(13x^2-10xy+13y^2-36x+36y)=0$

\begin{figure}[h] 
\begin{center}
\includegraphics[scale=0.45,angle=-90]{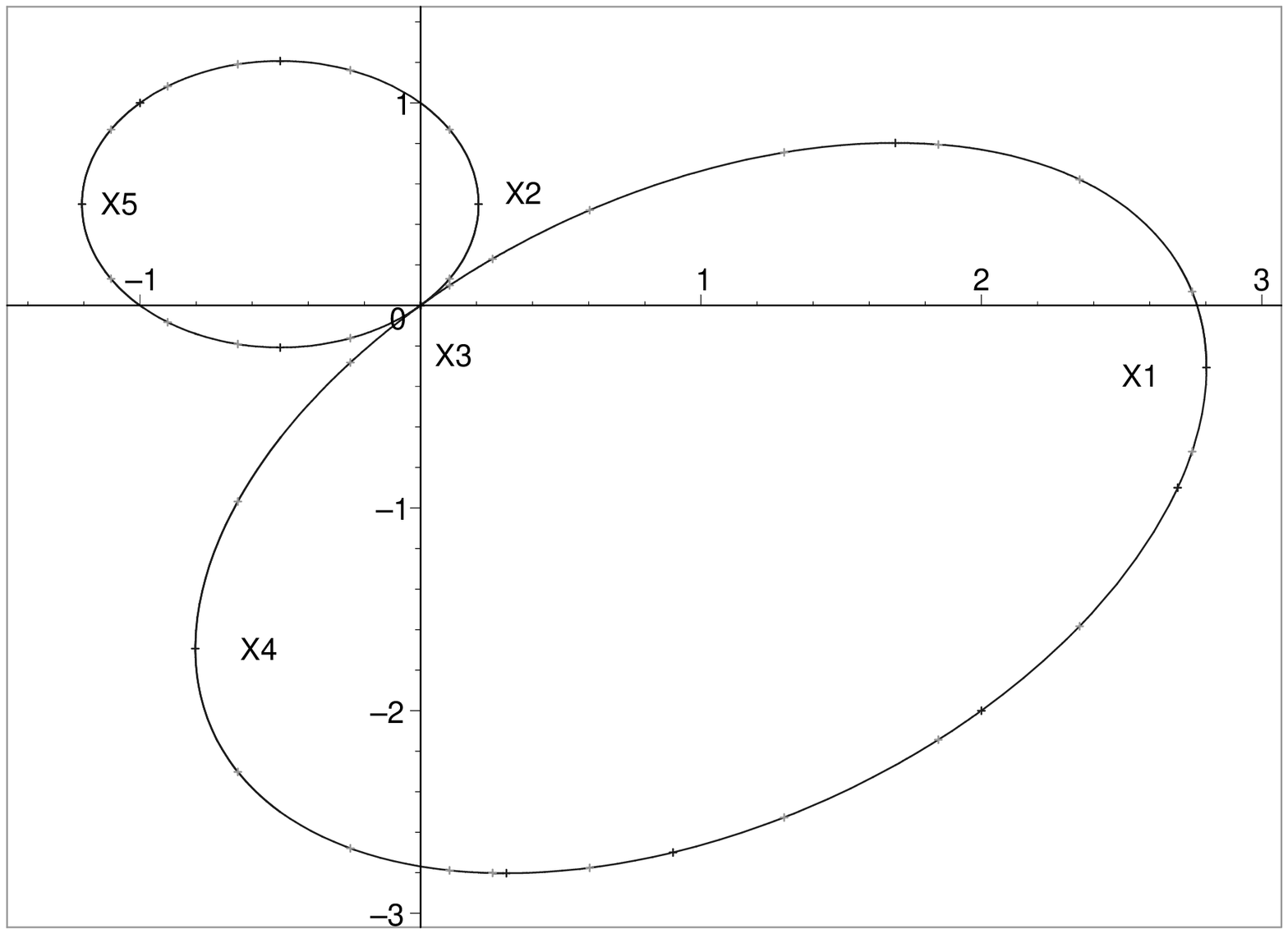}
\caption{Real part of $r(C)$ for Example $7$}
\label{fig:1m2a}
\end{center}
\end{figure}

\mytable{
\begin{center}
\begin{tabular}{||l||c|c|c|c||}
\hline
\hline
& & & & \\
Singular point & $x_1$ &   $x_2$ & $x_3$ & $x_4$\\ 
& & & & \\
\hline
& & & & \\
Braid monodromy &
$\s1^2$ & 
$\s2$ &
$\s3\s2\s1\s2^{-1}\s3^{-1}$ &
$\s3\s2^4\s3^{-1} $ \\
& & & & \\
\hline
%& & & & \\
\includegraphics[scale=1]{TitleGeometricHT.eps} &
\includegraphics[scale=0.5]{1x5.eps} &
\includegraphics[scale=0.5]{1x1.eps} &
\includegraphics[scale=0.5]{1x2.eps} &
\includegraphics[scale=0.5]{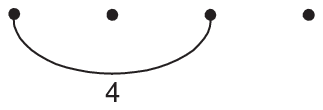} \\
\hline
\hline
& & & & \\
Singular point & $x_5$ & $x_6$ & $x_7$ & \\ 
& & & &\\
\cline{1-4}
& & & &\\
Braid monodromy &
$\s3^2\s2\s3^{-2}$ & 
$\s3\s2^{-1}\s1\s2\s3^{-1}$  & 
$\s3^2$ & \\
& & & & \\
\cline{1-4}
%& & & & \\
\includegraphics[scale=1]{TitleGeometricHT.eps} &
\includegraphics[scale=0.5]{7bX6.eps} &
\includegraphics[scale=0.5]{7bX3.eps} & 
\includegraphics[scale=0.5]{1x4.eps} &
\\
\hline
\hline
\end{tabular}
\caption{Braid monodromy results for Example $7$}\label{tab:1m2-a}
\end{center}}
\newpage

%%%%%%%%%%%%%%%%%%%%%%%%%%%%%%%%%%%%%%%%%%%%%%%%%%%%%%%%%%%%%%%%%%%%%%%%%%%%%%%%%%%%%%%%%%%%%%%
%                                                                                             %
% One intersection points of multiplicity $2$ - type b                                        %
%                                                                                             %
%%%%%%%%%%%%%%%%%%%%%%%%%%%%%%%%%%%%%%%%%%%%%%%%%%%%%%%%%%%%%%%%%%%%%%%%%%%%%%%%%%%%%%%%%%%%%%%
\subsection{Example 8. One intersection points of multiplicity 2 - type b} $ $\\
$p_0(x,y)=1+x^2+y^2$, $p_1(x,y)=4x+y+2x^2+y^2$, $p_2(x,y)=4x+y-2x^2-y^2$ \\
$r(C)$ is defined by: $(x^2+y^2-x+y)(5x^2-6xy+5y^2-16x+16y)=0$

\begin{figure}[h] 
\begin{center}
\includegraphics[scale=0.45,angle=-90]{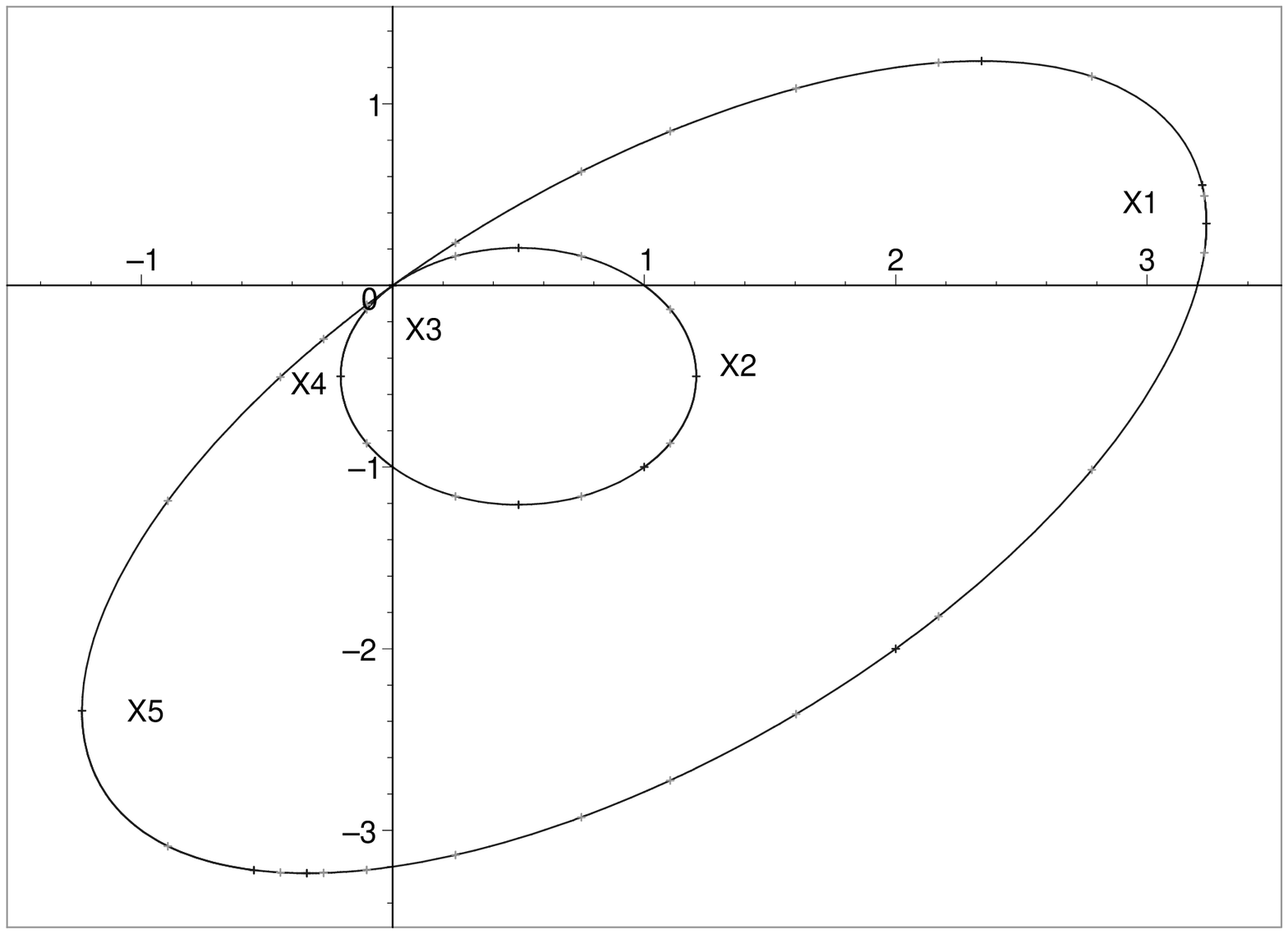}
\caption{Real part of $r(C)$ for Example $8$}
\label{fig:1m2b}
\end{center}
\end{figure}

\mytable{
\begin{center}
\begin{tabular}{||l||c|c|c|c||}
\hline
\hline
& & & & \\
Singular point & $x_1$ &  $x_2$ & $x_3$ & $x_4$\\ 
& & & & \\
\hline
& & & & \\
Braid monodromy &
$\s2$ &
$\s2^{-1}\s1^2\s2$ &
$\s3\s2^{-1}\s1\s2\s3^{-1}$ &
$\s1^4$ \\
& & & & \\
\hline
%& & & & \\
\includegraphics[scale=1]{TitleGeometricHT.eps} &
\includegraphics[scale=0.5]{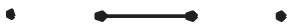} &
\includegraphics[scale=0.5]{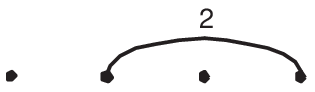} &
\includegraphics[scale=0.5]{7bX3.eps} & 
\includegraphics[scale=0.5]{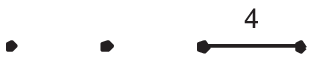} \\
\hline
\hline
& & & & \\
Singular point & $x_5$ & $x_6$ & $x_7$ & \\ 
& & & &\\
\cline{1-4}
& & & &\\
Braid monodromy &
$\s3\s2\s1\s2^{-1}\s3^{-1}$ & 
$\s3^2\s2\s3^{-2}$  & 
$\s3^2$ & \\
& & & & \\
\cline{1-4}
%& & & & \\
\includegraphics[scale=1]{TitleGeometricHT.eps} &
\includegraphics[scale=0.5]{1x2.eps} &
\includegraphics[scale=0.5]{7bX6.eps} & 
\includegraphics[scale=0.5]{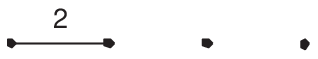} & \\
\hline
\hline
\end{tabular}
\caption{Braid monodromy results for Example $8$}\label{tab:1m2-b}
\end{center}}

%%%%%%%%%%%%%%%%%%%%%%%%%%%%%%%%%%%%%%%%%%%%%%%%%%%%%%%%%%%%%%%%%%%%%%%%%%%%%%%


\clearpage

\begin{thebibliography}{}

\bibitem{Artin} Artin, E., {\it Theory of braids}, 
Ann. Math. {\bf 48} (1947), 101-126.

\bibitem{BMGprept}
Ball, J. A., Malekorn, T. and Groenewalde, G. {\it Structured
noncommutative multidimensional linear systems}, preprint.

\bibitem{BR84}
Berstel, J. and Reutnaur, C. {\it Rational Series and their
Languages}, EATCS Monographs on Theoretical Computer Science,
Springer, 1984.

\bibitem{Birman} Birman, J., {\it Braids, links and mapping class groups}, 
Ann. Math Studies 82, Princeton University Press, 1975.

\bibitem{Dehornoy} Dehornoy, P., {\it Braids and Self Distributivity}, 
Progress in Mathematics, volume 192;. Birkhauser (2000).

\bibitem{F74}
Fliess, M., {\it Matrices de Hankel}, J. Math Pure Appl., 53,
(1974) and 197-222 \& erratum 54 (1975).

\bibitem{HMcCV}
Helton, W., McCullough, S. and Vinnikov, V., {\it Noncommutative
Convexity Arises from Linear Matrix Inequalities}, in preparation.

\bibitem{GeneralMT}
Kaplan, S. and Teicher, M., {\it Computing Braid Monodromy and the Moishezon Teicher Algorithm}, in preparation.

\bibitem{K80}
Kravitsky N., {\it On the discriminant function of two commuting
nonselfadjoint operators}, Integral Equations Operator Theory 3/1,
p. 97-124, 1980.

\bibitem{KuTe} Kulikov,~V.~S. and Teicher, M., {\it Braid monodromy
    factorizations and diffeomorphism types},
    Izv. Ross. Akad. Nauk Ser. Mat. {\bf 64}(2), 89--120 (2000) [Russian];
    {\it English transl.}, Izvestiya Math. {\bf 64}(2), 311--341 (2000).

\bibitem{BGTI} Moishezon, B. and Teicher, M.,
{\it Braid group techniques in complex geometry I, Line
arrangements in $\C \PP^2$}, Contemporary Math. 78 (1988),
425-555.

\bibitem{BGTII} Moishezon, B. and Teicher, M.,
{\it Braid Group Techniques in Complex Geometry II: From
arrangements of lines to cuspidal curves}, {\it LNM} 1479 (1989),
131-179.

\bibitem{S61}
Sch\"{u}tzenberger, M. P. {\it On the definition of a family of
automata}, Information and Control 4 1961 p. 245--270.

\bibitem{SVprept}
Shapiro A., Vinnikov V., {\it Rational transformations of
algebraic curves and elimination theory}, Linear Algebra and its
Applications, preprint.

\bibitem{V93}
Vinnikov V., {\it Self-adjoint determinantal representations of
real plane curves}, Math. Ann., 296 (1993), p. 453-473.

\end{thebibliography}
\end{document}